\documentclass{conm-p-l}

\copyrightinfo{2011}{}

\usepackage{graphicx}
\usepackage{subfigure}

\usepackage{amssymb,amsmath,amsthm,amscd}

\theoremstyle{definition}

\theoremstyle{remark}

\numberwithin{equation}{section}

\newcommand{\tu}{\tilde{u}}
\newcommand{\tU}{\tilde{U}}

\newcommand{\e}{\epsilon}

\newcommand{\Dl}{\Delta}

\newcommand{\ra}{\rightarrow}

\newcommand{\sg}{\sigma}

\newcommand{\pa}{\partial}

\newcommand{\la}{\lambda}

\newcommand{\Om}{\Omega}

\begin{document}

\title[A resolution of the turbulence paradox]
{A resolution of the turbulence paradox: numerical implementation}

\author{Yueheng Lan}
\address{Department of Physics, Tsinghua University, Beijing 100084, China}
\email{yueheng\_lan@yahoo.com}

\author{Y. Charles Li}
\address{Department of Mathematics, University of Missouri, 
Columbia, MO 65211, USA}
\email{liyan@missouri.edu}

\curraddr{}
\thanks{}

\date{}

\dedicatory{}

\keywords{Sommerfeld paradox, Couette flow, flow stability, shear flow, transition to turbulence.}

\begin{abstract}
Sommerfeld paradox (turbulence paradox) roughly says that mathematically the Couette linear shear flow 
is linearly stable for all values of the Reynolds number,
but experimentally transition from the linear shear to turbulence occurs under perturbations of any size 
when the Reynolds number is large enough. In \cite{LL09}, we offered a resolution of this paradox. The 
aim of this paper is to provide a numerical implementation of the resolution. The main idea of the 
resolution is that even though the linear shear is linearly stable, slow orbits (also called quasi-steady 
states) in arbitrarily small neighborhoods of the linear shear can be linearly unstable. The key is that 
in infinite dimensions, smallness in one norm does not mean smallness in all norms.
Our study here focuses upon
a sequence of 2D oscillatory shears which are the Couette linear shear plus small amplitude and high 
frequency sinusoidal shear perturbations. In the fluid velocity 
variable, the sequence approaches the Couette linear shear (e.g. in $L^2$ norm of velocity), thus it can 
be viewed as Couette linear shear plus small noises; while in the fluid vorticity variable, the sequence does 
not approaches the Couette linear shear (e.g. in $H^1$ norm of velocity). Unlike the Couette linear shear, 
the sequence of oscillatory shears has inviscid linear instability; 
furthermore, with the sequence of oscillatory shears as potentials, the Orr-Sommerfeld operator has 
unstable eigenvalues when the 
Reynolds number is large enough, this should lead to transient nonlinear growth which manifests as 
transient turbulence as observed in experiments. The main result of this paper verifies this transient growth. 
\end{abstract}

\maketitle

\section{Introduction}

The most influential paradox in fluids is the d'Alembert paradox saying that a body moving through 
water has no drag as calculated by d'Alembert \cite{dAl52} via inviscid theory, while experiments 
show that there is a substantial drag on the body. The paradox splitted the field of fluids into 
two branches: 1. Hydraulics --- observing phenomena without mathematical explanation, 2. Theoretical 
Fluid Mechanics --- mathematically predicting phenomena that could not be observed. A revolutionary 
development of the boundary layer theory by Ludwig Prandtl in 1904 resolved the paradox by paying 
attention to the substantial effect of small viscosity in the boundary layer \cite{Pra04}. Prandtl's 
boundary layer theory laid the foundation of modern unified fluid mechanics. 

Sommerfeld paradox has the potential of being the next most influential paradox in fluids. The paradox says 
that the Couette linear shear flow is linearly stable for all values of the Reynolds numbers as first 
calculated by Sommerfeld \cite{Som08}, 
but experiments show that transition from the linear shear to turbulence occurs under perturbations of any size 
when the Reynolds number is large enough \cite{Bec95} \cite{Kre94}. This paradox 
lies at the heart of understanding turbulence inside the infinite dimensional phase space. Dynamical 
system studies on the Navier-Stokes flow in an infinite dimensional phase space is still at its developing 
stage \cite{LL08}. A typical dynamical system study on chaos often starts from fixed points (steady states) and then
pursues their invariant manifolds to understand the phase space structures. Certain techniques e.g. Melnikov 
integrals can then be used to detect the intersection between invariant manifolds and the existence of 
homoclinic or heteroclinic orbits. In some cases, chaos (turbulence) can be rigorously proved to exist 
using e.g. shadowing techniques \cite{Li04}. For the Navier-Stokes equations, the problem is much more difficult 
than a typical dynamical system. The Sommerfeld paradox is a clear example.

The status of research on transition to turbulence may be stated as follows: (1). The cause (initiator) of the 
transition should have general principles and can be clarified. (2). After the initiation of the transition, 
further development of turbulence may not have any general principle. Such further developments may differ 
substantially between 2D and 3D turbulence. Due to the extra constant of motion given by the enstrophy, invariant 
manifolds of 2D Euler equations, if they exist, may be degenerate, i.e. stable and unstable manifolds coincide, 
so that 2D turbulence develops mildly and slowly, while 3D turbulence develops violently and quickly.

Besides the simple steady states like the linear shear in Couette flow, there are also numerical pursuits on 
more general steady states. For plane Couette flow, plane Poiseuille flow, and pipe Poiseuille flow, unstable 
steady states with three dimensional spatial patterns have been intensively explored numerically \cite{Nag90} 
\cite{OK80} \cite{OP80} \cite{Eck08} \cite{Ker05}  \cite{Vis08}. Periodic orbits \cite{KK01} and quasi-periodic orbits 
\cite{Vis07} in plane Couette flow are also explored. 

Explorations on two dimensional non-trivial steady states turn out to be not successful \cite{ENR08}. That is, 
the counterpart of the 3D upper or lower branch steady state \cite{Nag90} has not been found in 2D. On the other hand, 
numerics shows that transitions still occur from the linear shear to turbulence in 2D.

There have been a lot of studies on the problem of transition to turbulence \cite{SH01}. There are also some previous 
attempts to explain the Sommerfeld paradox. One popular attempt which was first suggested by Orr \cite{Orr07}, is to use 
the non-normality of the linearized
Navier-Stokes operator to get algebraic growth of perturbations before their
final decay. (Note: non-normality refers to operators with non-orthogonal
eigenfunctions.) However, it is not clear how such linear algebraic growth
relates to the nonlinear dynamics. Moreover, the
non-normality theory cannot explain many coherent structures observed in the
transient turbulence.

Here we mainly focus upon a sequence of 2D oscillatory shears in plane Couette flow. These oscillatory shears 
are built from single Fourier 
mode modifications to the linear shear. The sequence of 2D oscillatory shears approaches the linear shear 
in the velocity variable but not the vorticity variable. All these oscillatory shears have a linear 
instability. We believe that such an 
instability explains the initiation of transition from the linear shear to turbulence. The main numerical result of 
this paper verifies this claim. Near the 2D oscillatory shears, inviscid 
two dimensional steady states (with a cat's eye structure) can be established rigorously \cite{LL09}. Coherent structures 
revealed in the current study are slightly viscous continuations of the inviscid cat's eye steady states.

Of course, our oscillatory shears are linearly unstable to 3D perturbation modes too. In \cite{LL09}, 
we show that 3D shears in a neighborhood of our oscillatory shears are linearly unstable too. 
The sequence of 2D oscillatory shears focused upon here is just a representative of all linearly unstable slow orbits 
in the neighborhood of the linear shear. There are many other such linearly unstable slow orbits in the neighborhood 
of the linear shear \cite{Li11a} \cite{Li11b}. All these linear instabilities contribute to the initiation of the 
transition fom the linear shear to turbulence!

\section{The sequence of oscillatory shears: analytical results}

Two dimensional plane Couette flow is governed by the Navier-Stokes equations with specific boundary conditions,
\begin{equation}
u_{i,t} + u_j u_{i,j} = -p_{,i} +\e u_{i,jj} \ , \quad u_{i,i} = 0 \ ; 
\label{NS}
\end{equation}
defined in a horizontal channel, where $u_i$ ($i=1,2$) are the velocity 
components along $x$ and $y$ directions, $p$ is the pressure, and $\e$ is the inverse of the Reynolds number 
$\e = 1/R$; with the
boundary condition
\begin{equation}
u_1(t, x, 0) = 0 , \ u_1(t, x, 1) = 1 , \ u_2(t, x, 0) = u_2(t, x, 1) = 0;
\label{BC}
\end{equation}
and all the variables are periodic in $x$. The Couette linear shear 
($u_1 = y$, $u_2=0$) is linearly and nonlinearly stable for 
any Reynolds number \cite{Rom73} \cite{HKK09}. Here we focus on the sequence of oscillatory shears
\begin{equation}
u_1 = U(y) = y + \frac{c}{n} \sin (4n\pi y),  \quad  u_2=0 ; 
\label{Os}
\end{equation}
where $c$ is a constant. One can view this as a single mode of the Fourier series of all 2D shears satisfying 
the above boundary conditions (\ref{BC}),
\[
y + \sum_{m=1}^{+\infty} c_m \sin (my) .
\]
As $n \ra \infty$, the oscillatory shears approach the linear shear $U(y) \ra y$. On the 
other hand, 
in the vorticity variable, the oscillatory shears do not approach the linear shear $U_y = 1 + 4c\pi 
\cos (4n\pi y) \not\ra 1$.
Thus in the velocity variables, the oscillatory shears can be viewed as the linear shear plus small noises.
It is shown in \cite{LL09} that these oscillatory shears are linearly unstable under the 2D Euler flow when 
\begin{equation}
\frac{1}{2} \frac{1}{4\pi} < c < \frac{1}{4\pi}.
\label{cr}
\end{equation}
It is shown in \cite{LL09} that the oscillatory shears bifurcate into steady states of 2D Euler flow 
with cat's eye structures.
Under the 2D Navier-Stokes flow, these oscillatory shears are not steady, rather slowly drifting,
\begin{equation}
U(t,y) = y + \frac{c}{n} e^{-\e (4n\pi )^2 t}\sin (4n\pi y). 
\label{sdr}
\end{equation}
Nevertheless, it is shown in \cite{LL09} that  by simply using the  oscillatory shears  (\ref{Os}) under 
condition (\ref{cr}) as the potential,  the Orr-Sommerfeld operator (linear Navier-Stokes operator) has unstable 
eigenvalues which converge to those of linear Euler operator as $\e \ra 0^+$. Based upon the geometric intuition of 
the geometric singular perturbation theory of Fenichel \cite{Fen79} \cite{Fen74} \cite{Fen77} \cite{Fen71}, such a 
linear instability should lead to a transient nonlinear growth near the oscillatory shears 
(and the linear shear). The main goal of the current paper is to verify this numerically.
Such a growth will manifest itself in experiments and numerical simulations as a transition from the linear shear to 
turbulence. Here the amplitude 
of the perturbation from the linear shear will be measured by the deviation of the oscillatory shears from the linear shear and 
the perturbation on top of the oscillatory shears. Thus the main idea of our resolution of the Sommerfeld paradox 
\cite{LL09} is that even though the linear shear itself is linearly stable, slow orbits in arbitrarily small 
neighborhoods of the linear shear can be linearly unstable and lead to transitions from the linear shear to 
turbulence. The drifting shears (\ref{sdr}) are typical examples of such linearly unstable slow orbits. 

\section{The sequence of oscillatory shears: numerical results}

For the numerical simulations, we choose the x-direction spatial period $L_x = 2.2 \pi$.
We choose the initial conditions by adding random perturbations to the 
oscillatory shears (\ref{Os}). The $L^2$ norm of the random perturbations is $0.01$. Notice from 
(\ref{cr}) that the parameter $c$ is in the range ($0.04, 0.08$). For all the simulations, we choose 
$c=0.07$. In order for the random perturbations not to overwhelm the oscillatory part of our 
oscillatory shear (\ref{Os}), $n$ needs to be less than $7$. Here we shall simulate $n=1,2,3$. 
\begin{figure}[ht] 
\centering
\subfigure[$t=0$]{\includegraphics[width=2.3in,height=2.3in]{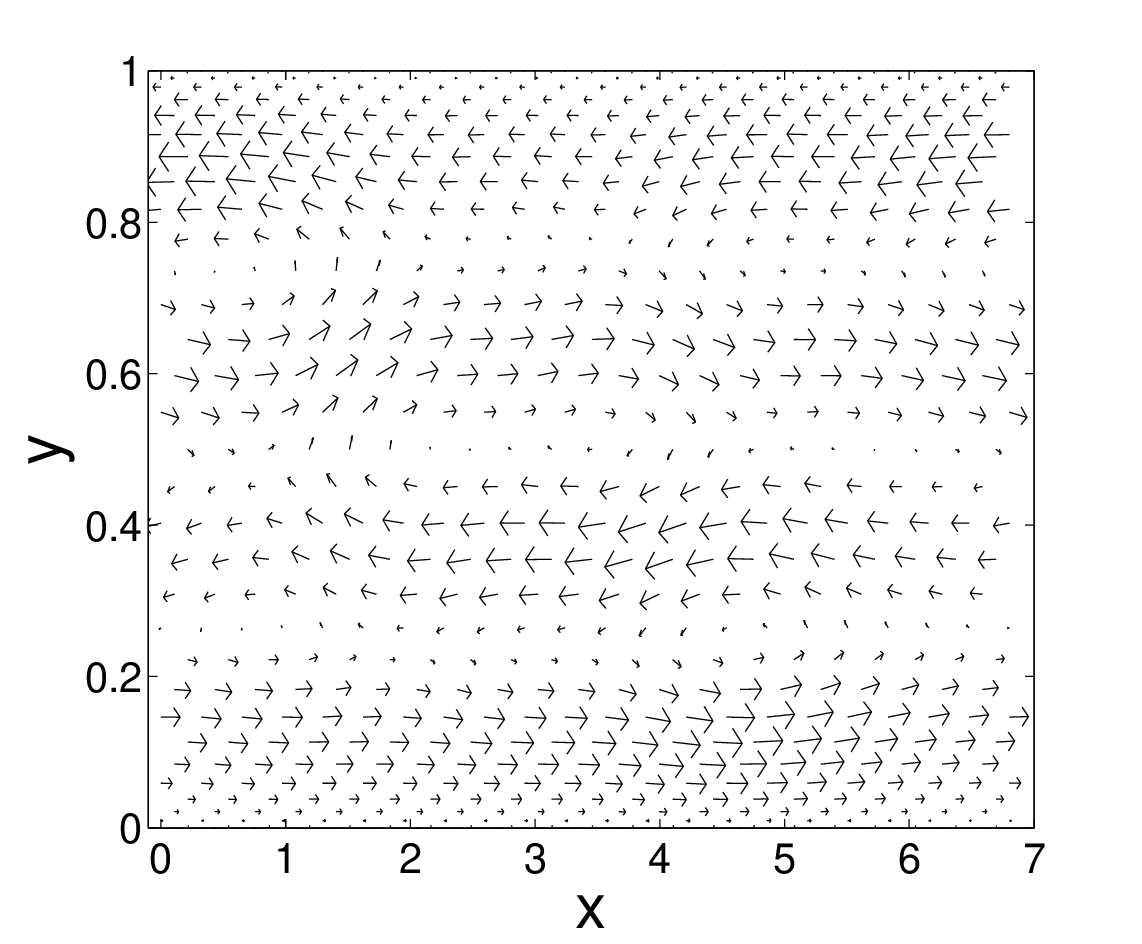}}
\subfigure[$t =4.9$]{\includegraphics[width=2.3in,height=2.3in]{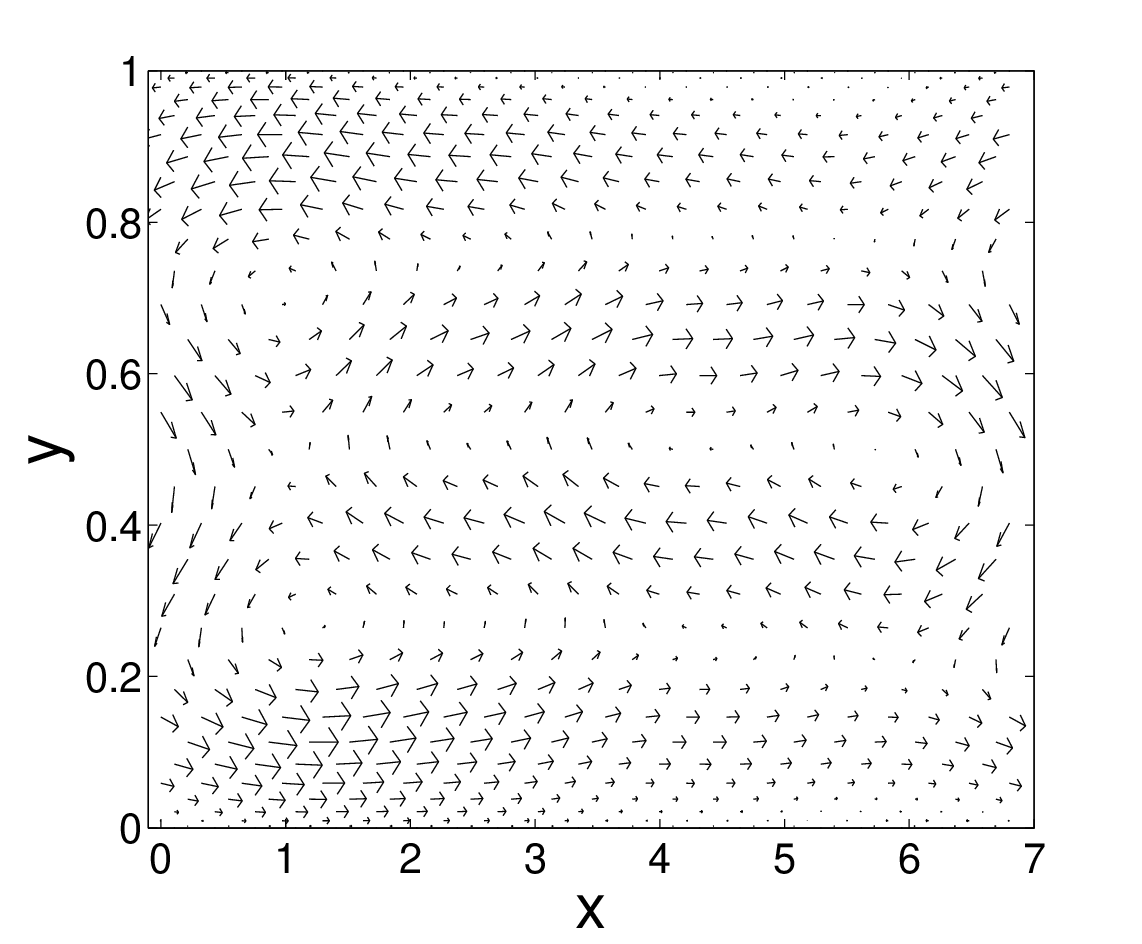}}
\subfigure[$t=14.9$]{\includegraphics[width=2.3in,height=2.3in]{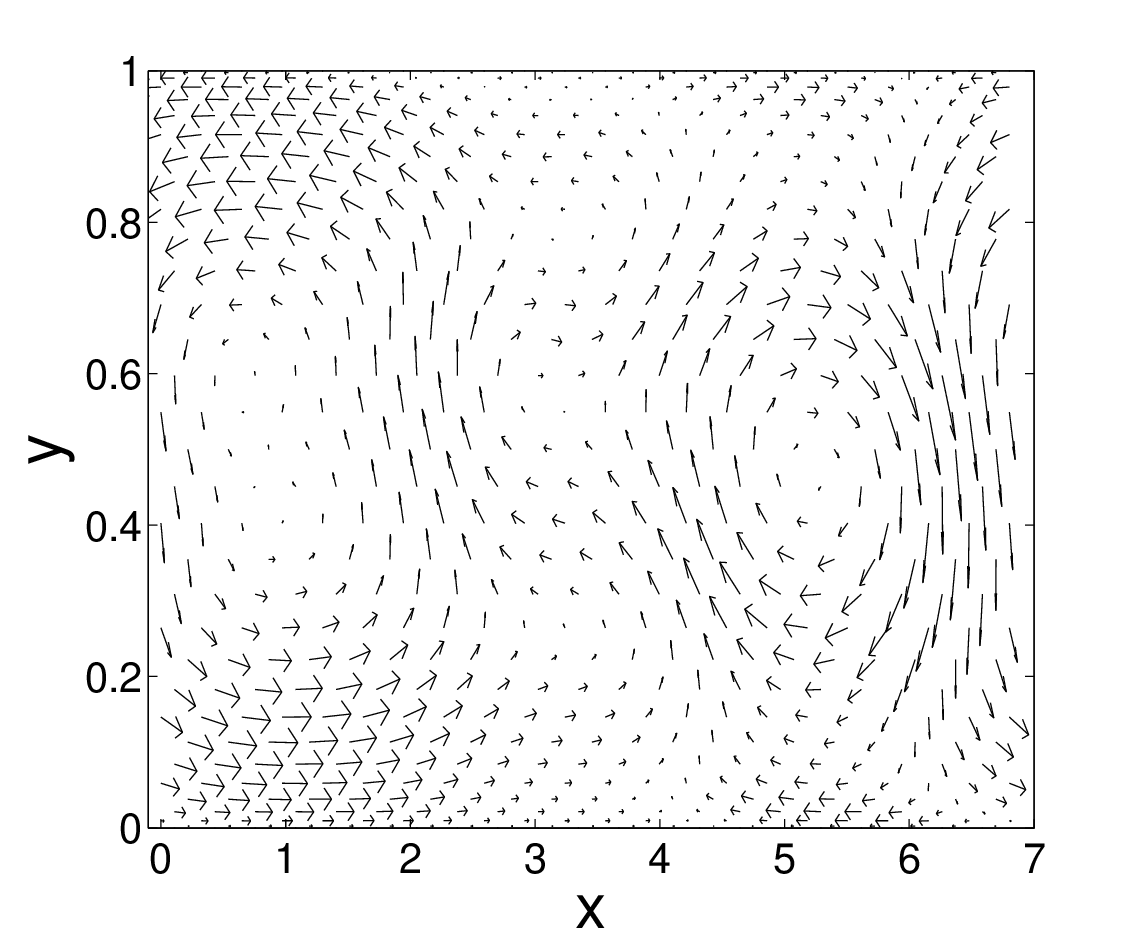}}
\subfigure[$t=24.9$]{\includegraphics[width=2.3in,height=2.3in]{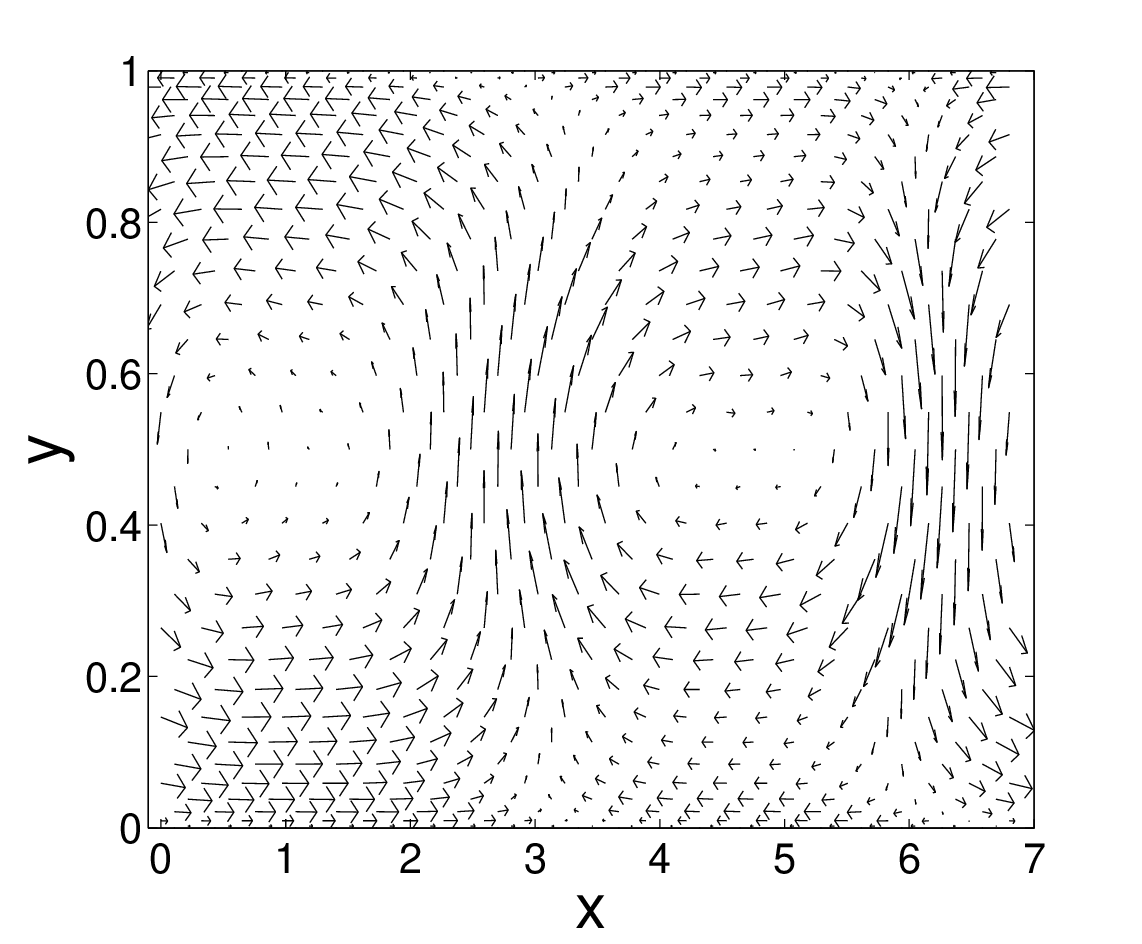}}
\caption{The development of a coherent structure and transient turbulence with initial condition (\ref{Os}) plus a random perturbation
(of $L^2$ norm $0.01$) for $n=1$, $c=0.07$, and $R=10000$. Here the linear shear has been subtracted.}
\label{cod}
\end{figure}
\begin{figure}[ht] 
\centering
\subfigure[$n=1, R=10000$]{\includegraphics[width=1.6in,height=1.6in]{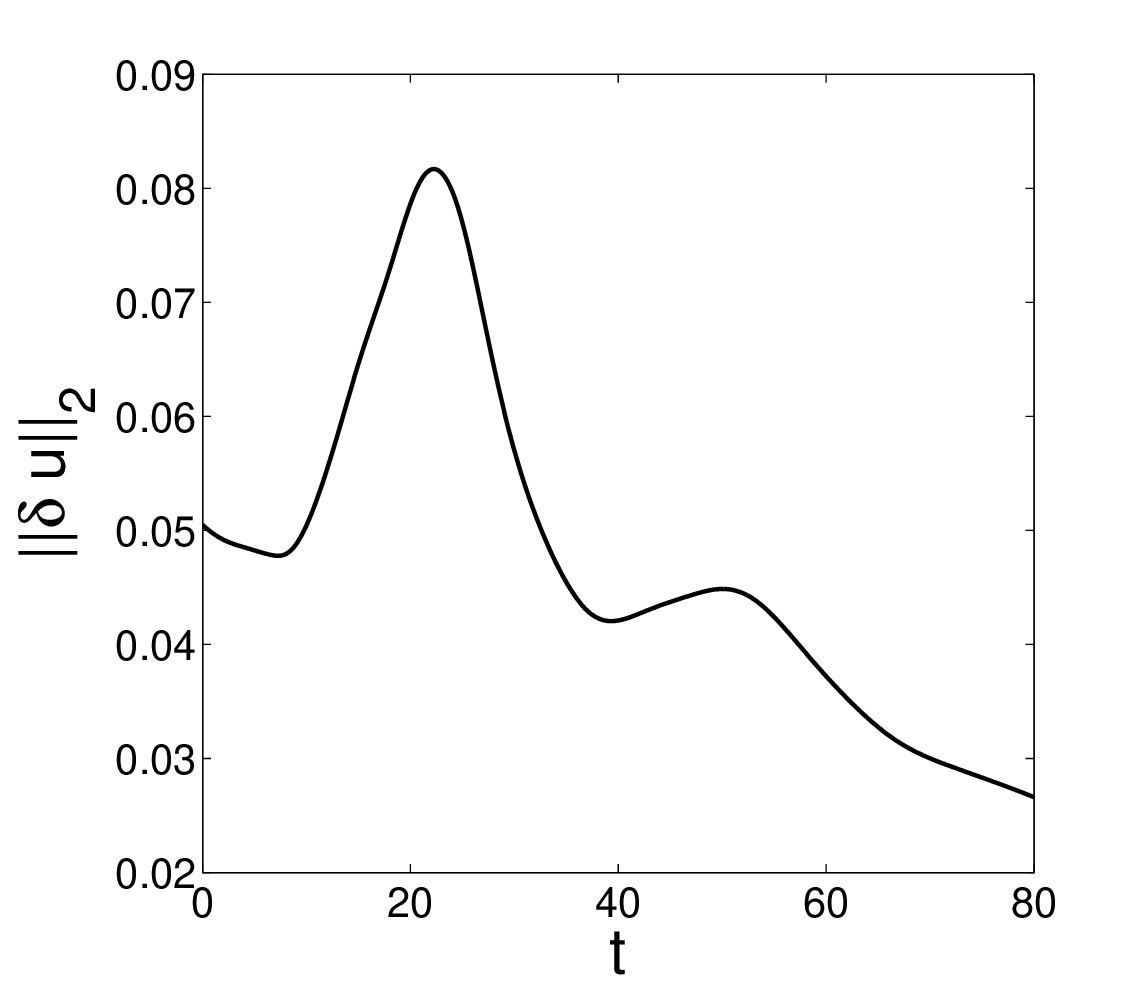}}
\subfigure[$n=1, R=20000$]{\includegraphics[width=1.6in,height=1.6in]{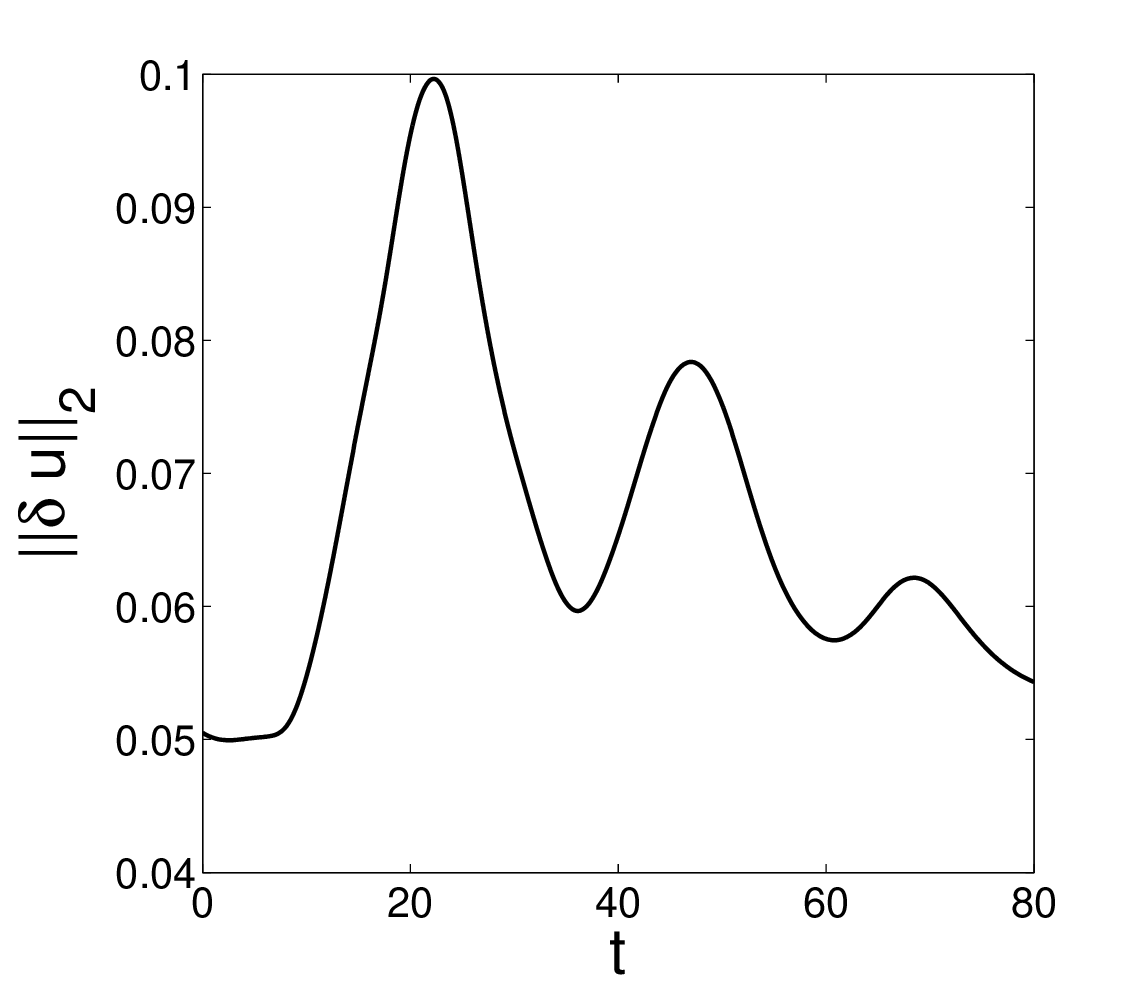}}
\subfigure[$n=2, R=10000$]{\includegraphics[width=1.6in,height=1.6in]{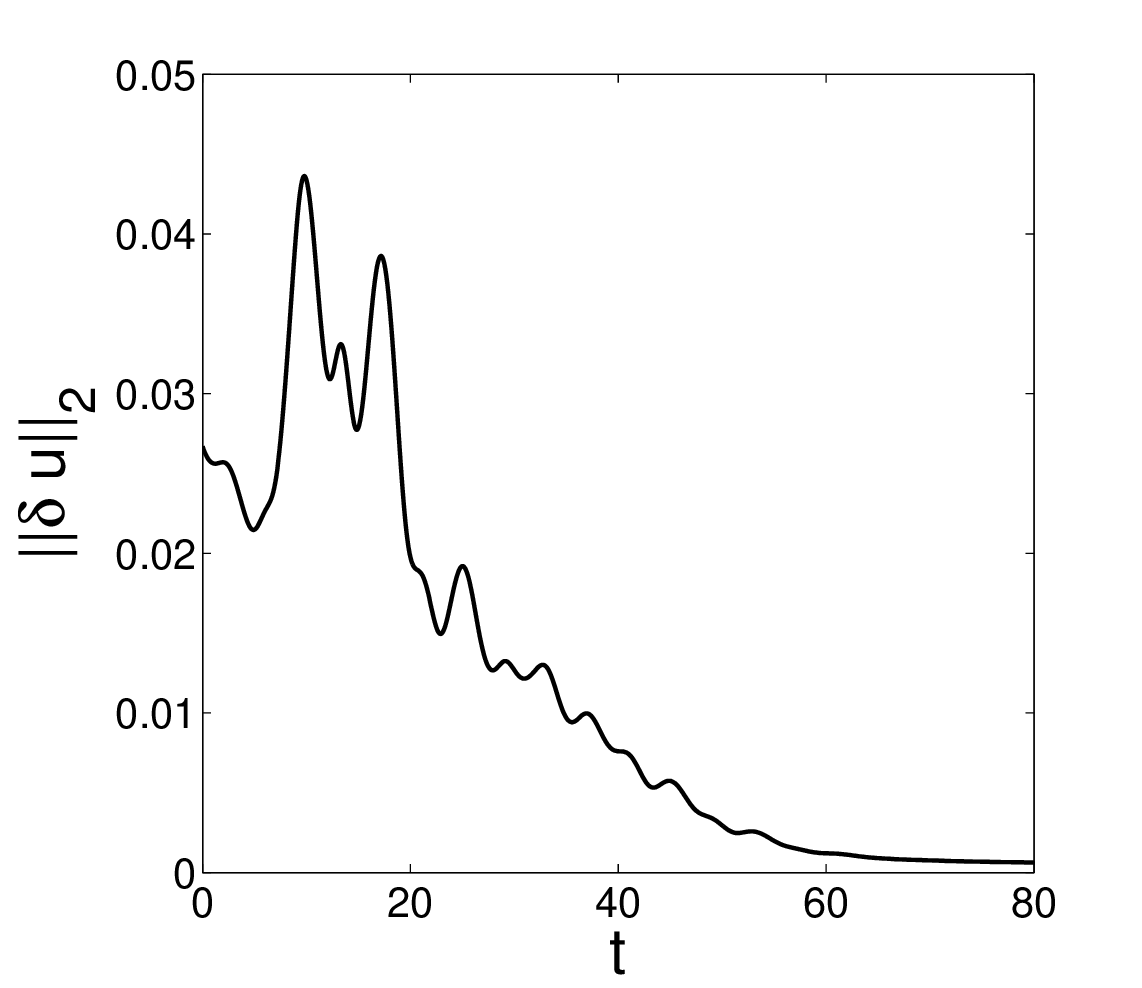}}
\subfigure[$n=2, R=20000$]{\includegraphics[width=1.6in,height=1.6in]{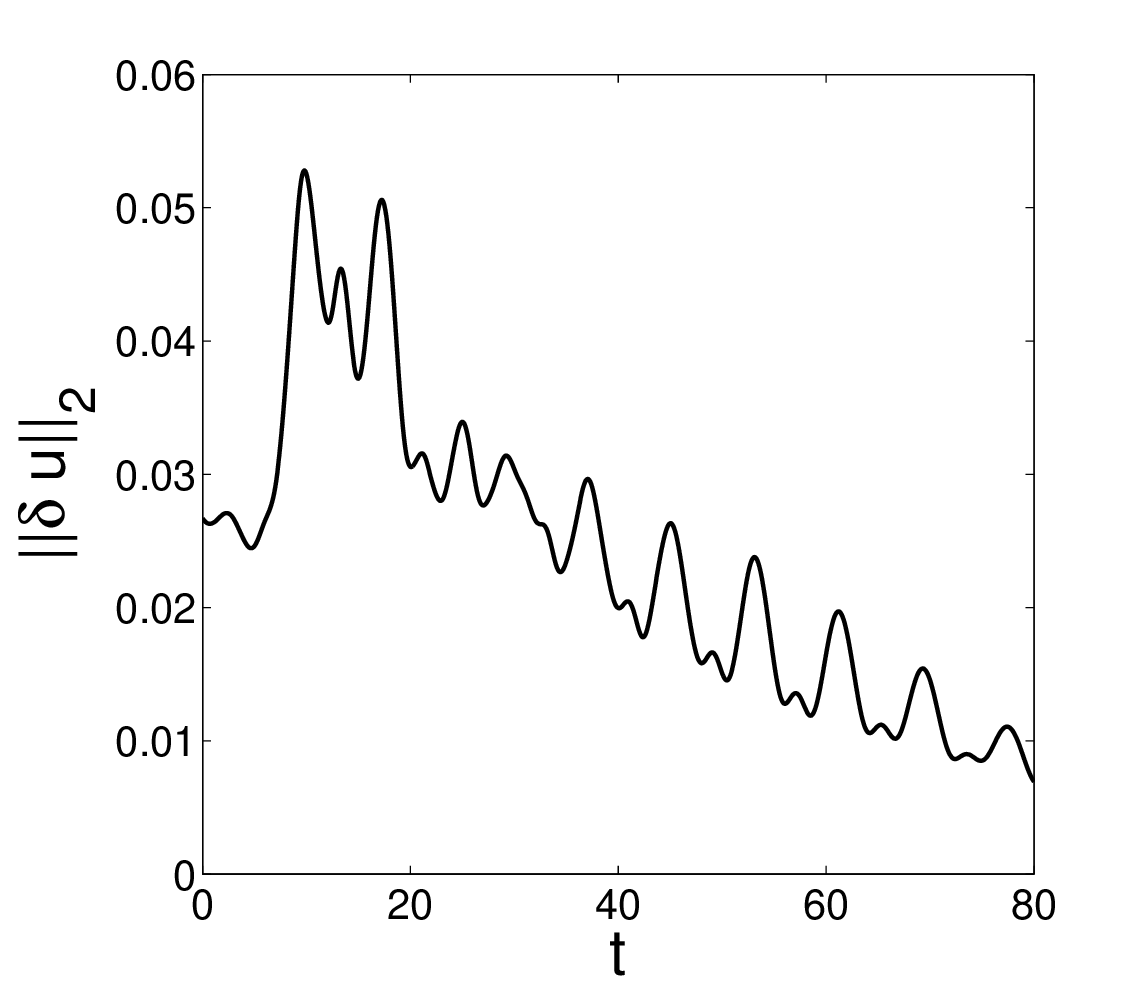}}
\subfigure[$n=3, R=10000$]{\includegraphics[width=1.6in,height=1.6in]{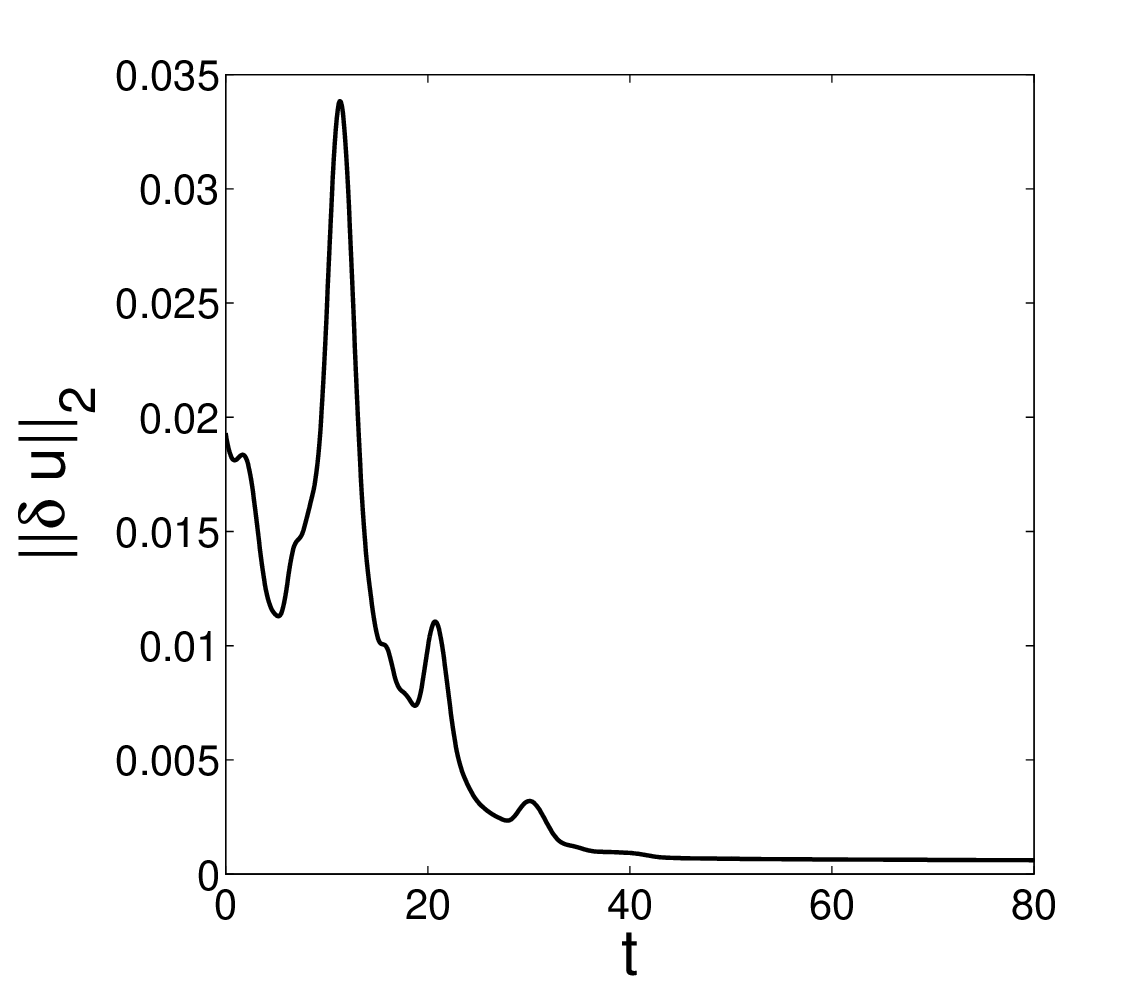}}
\subfigure[$n=3, R=20000$]{\includegraphics[width=1.6in,height=1.6in]{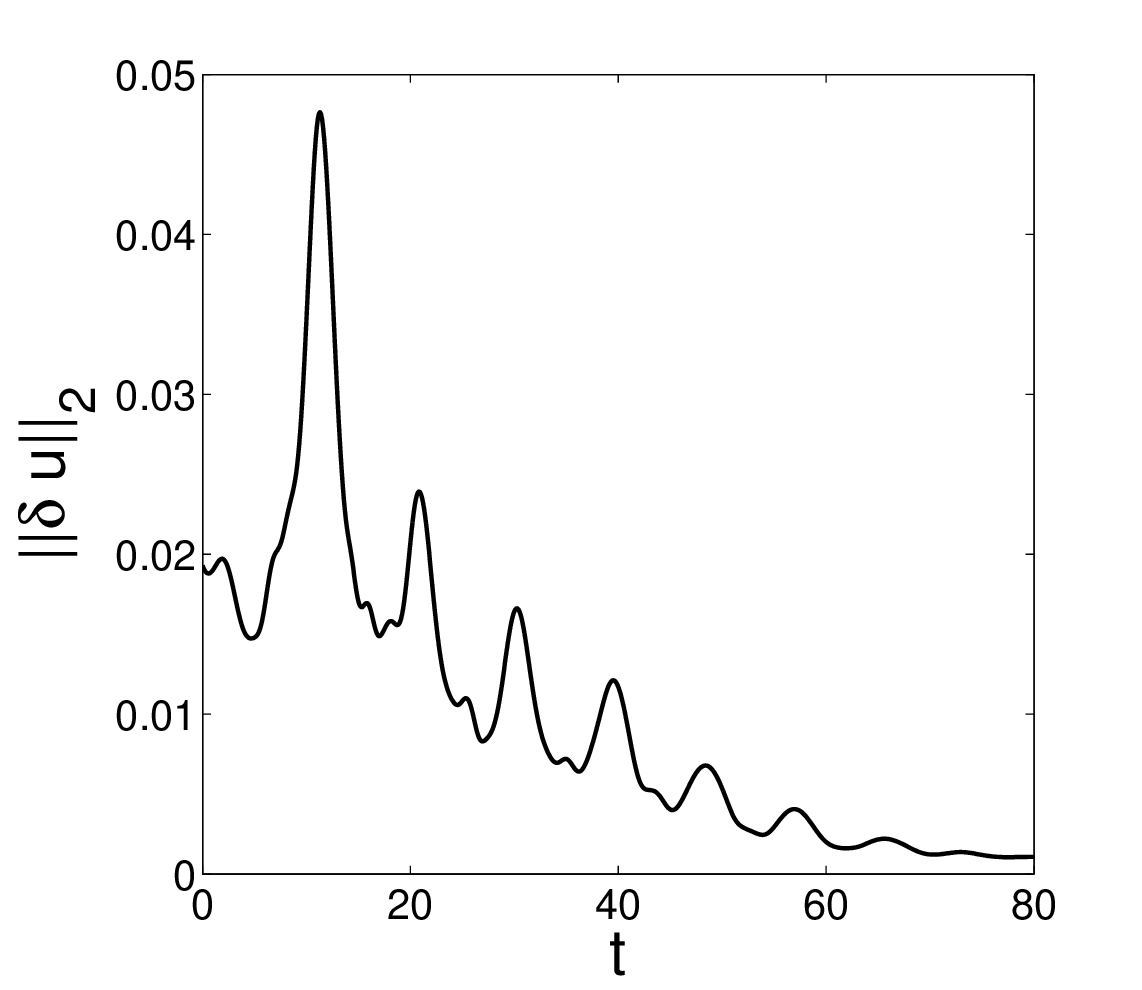}}
\caption{The growth of the $L^2$ norm of the deviation from the linear shear with 
initial condition (\ref{Os}) plus a random perturbation (of $L^2$ norm $0.01$). The growth rate $\sg$ is defined as the $\log_e$ of quotient (of the first maximum and the first minimum) divided by the time spent. (a). $\sg = 0.044$,
(b). $\sg = 0.055$, (c). $\sg = 0.24$, (d). $\sg = 0.25$, (e). $\sg = 0.17$, (f). $\sg = 0.18$.}
\label{grth}
\end{figure}
\begin{figure}[ht] 
\centering
\subfigure[$n=1, R=10000$]{\includegraphics[width=1.6in,height=1.6in]{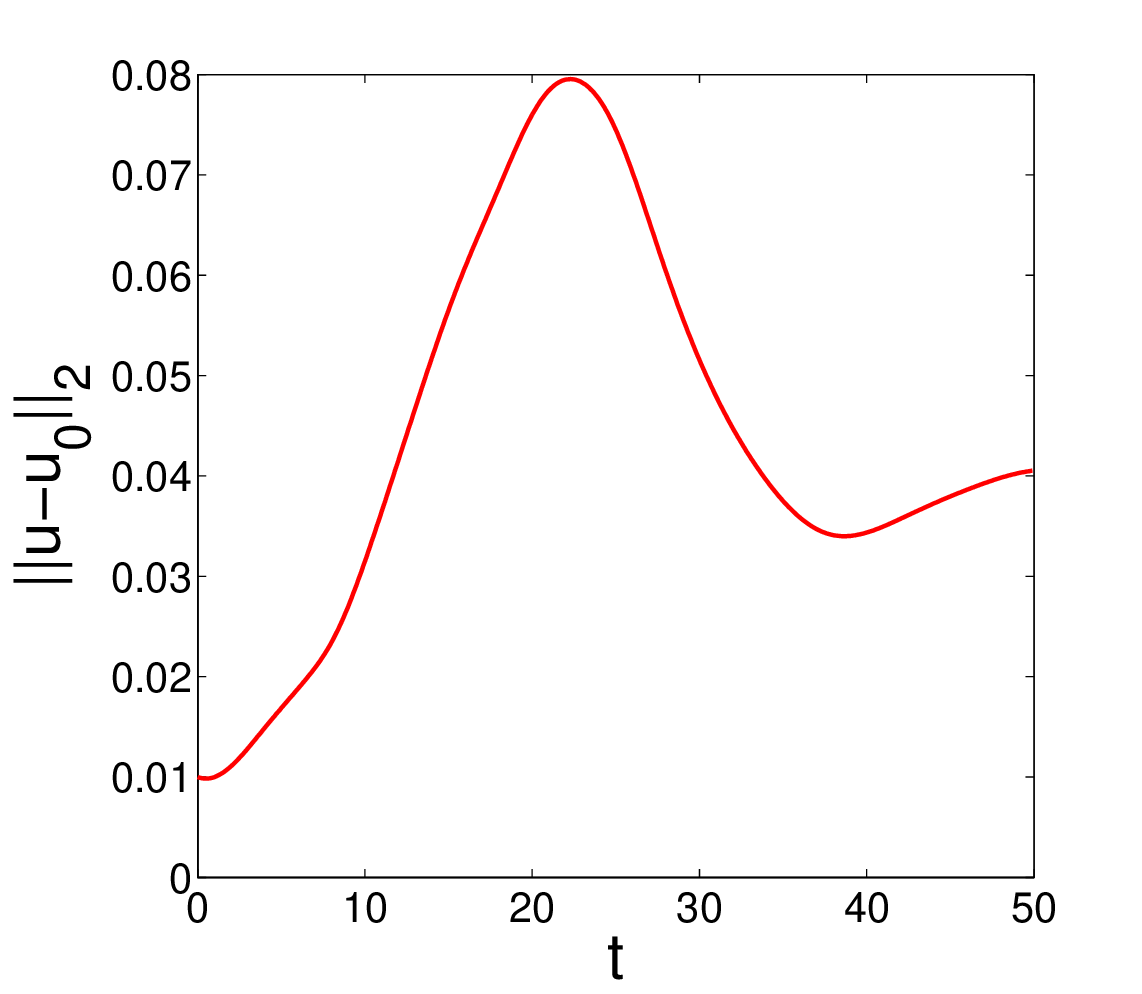}}
\subfigure[$n=1, R=20000$]{\includegraphics[width=1.6in,height=1.6in]{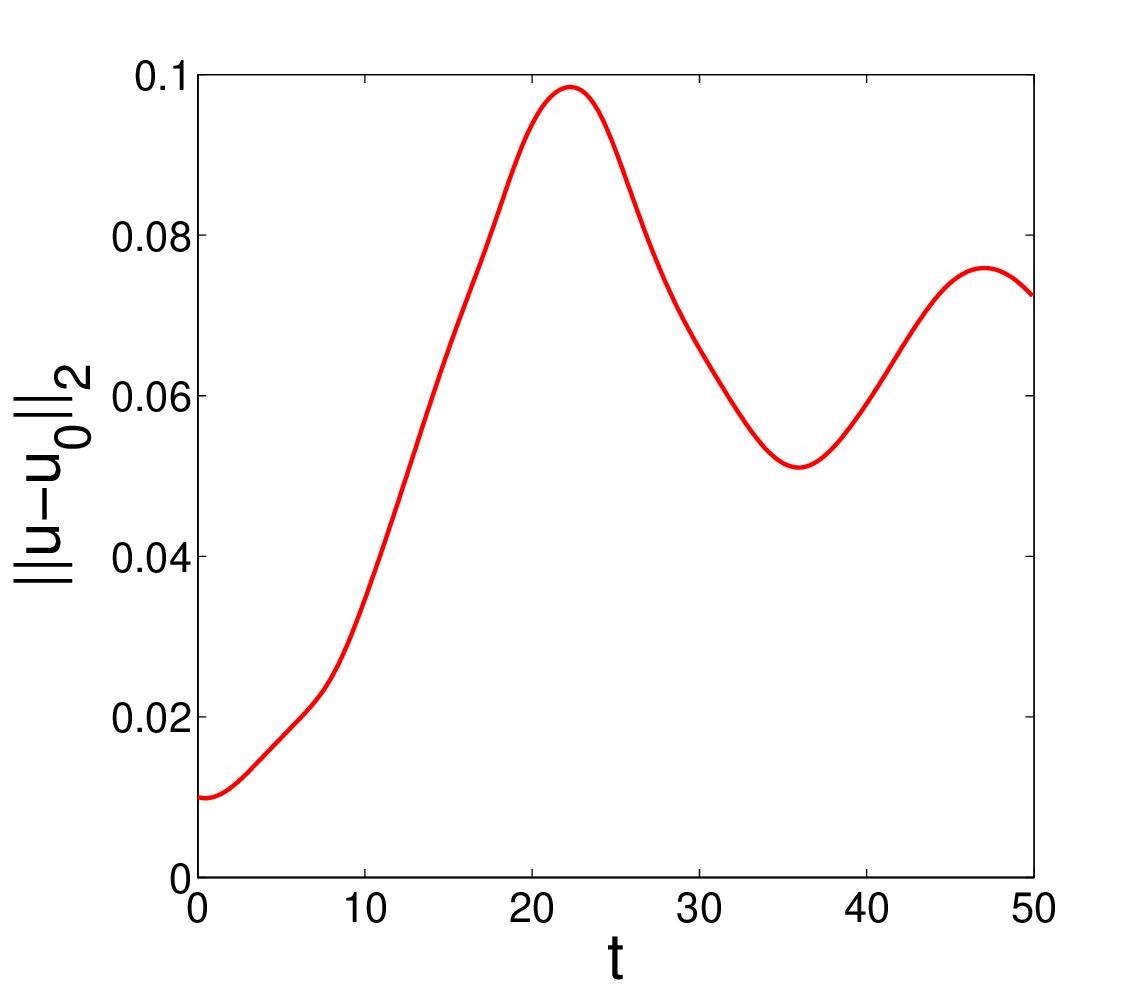}}
\subfigure[$n=2, R=10000$]{\includegraphics[width=1.6in,height=1.6in]{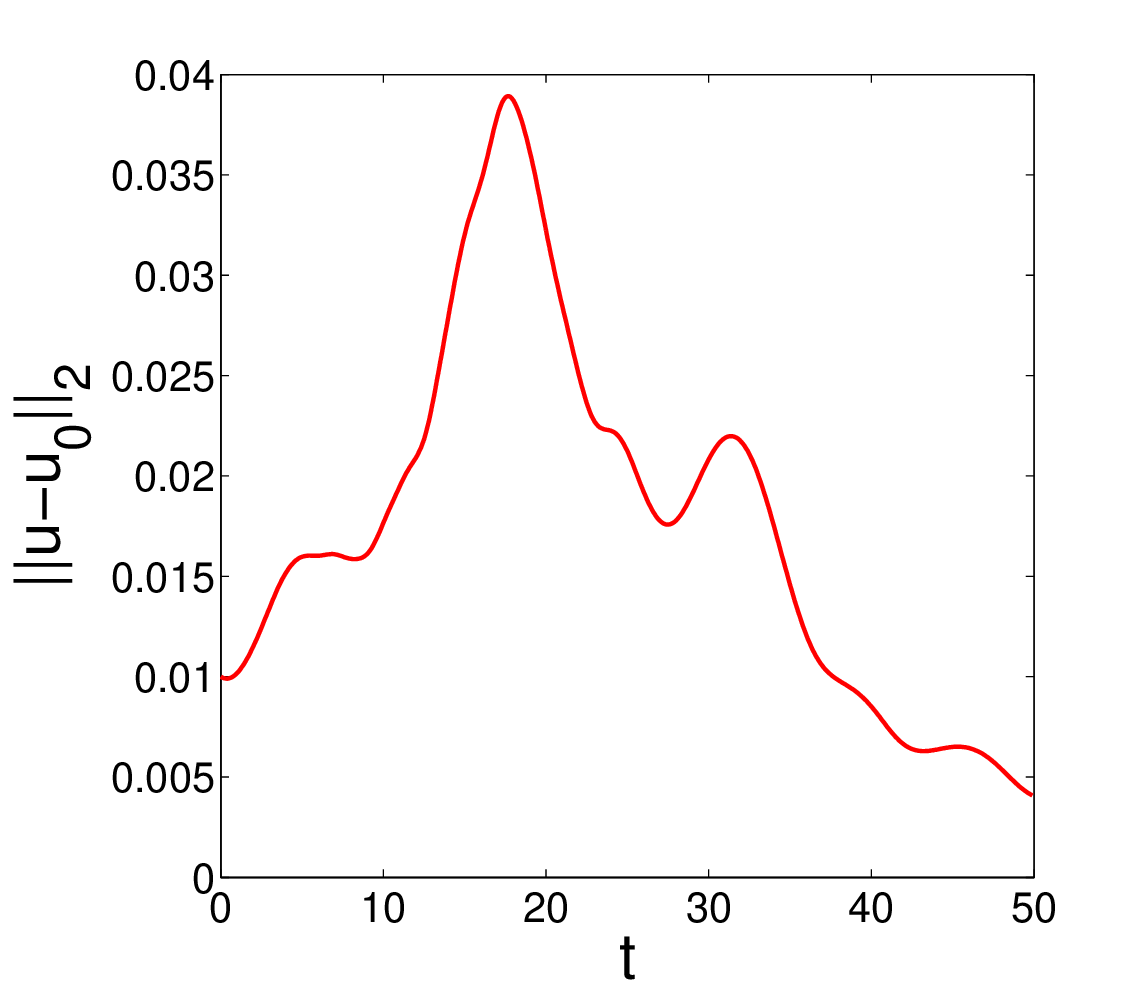}}
\subfigure[$n=2, R=20000$]{\includegraphics[width=1.6in,height=1.6in]{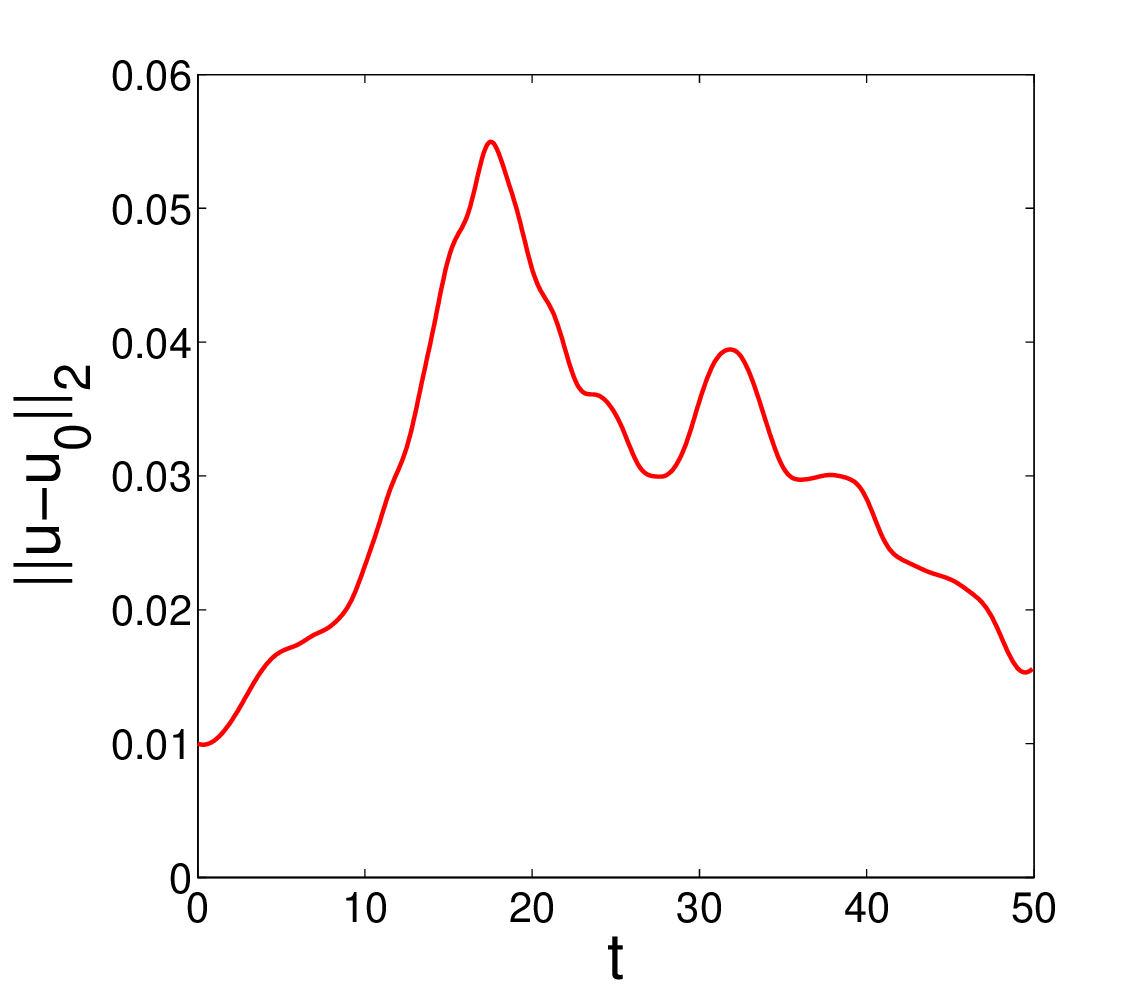}}
\subfigure[$n=3, R=10000$]{\includegraphics[width=1.6in,height=1.6in]{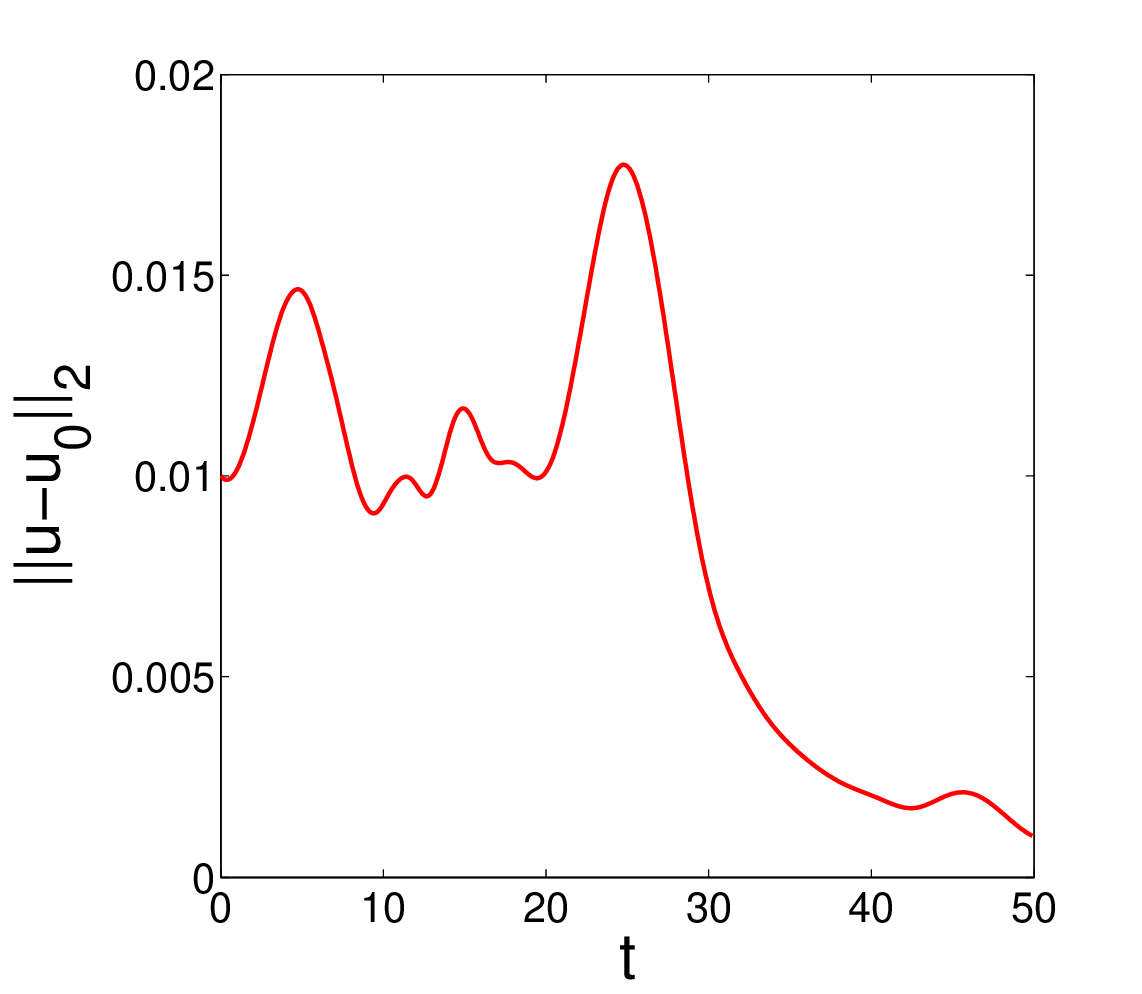}}
\subfigure[$n=3, R=20000$]{\includegraphics[width=1.6in,height=1.6in]{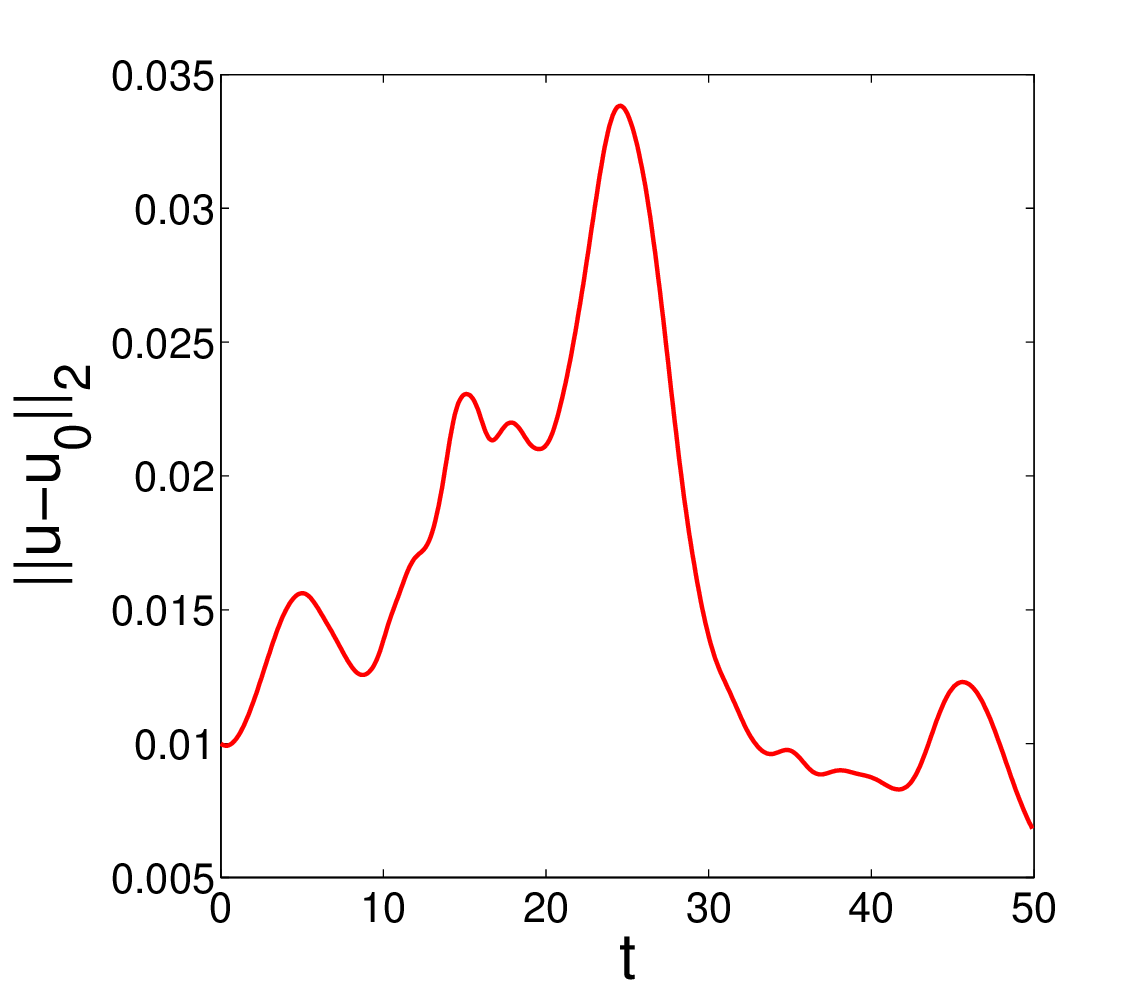}}
\caption{The growth of the $L^2$ norm of the deviation from the slow drifting $u_0$ (\ref{sdr}) with 
initial condition (\ref{Os}) plus a random perturbation (of $L^2$ norm $0.01$). The growth rate $\sg$ is defined as the $\log_e$ of quotient (of the first maximum and the first minimum) divided by the time spent. (a). $\sg = 0.11$,
(b). $\sg = 0.12$, (c). $\sg = 0.13$, (d). $\sg = 0.14$, (e). $\sg = 0.11$, (f). $\sg = 0.13$.}
\label{grth2}
\end{figure}
\begin{figure}[ht] 
\centering
\subfigure{\includegraphics[width=1.6in,height=1.6in]{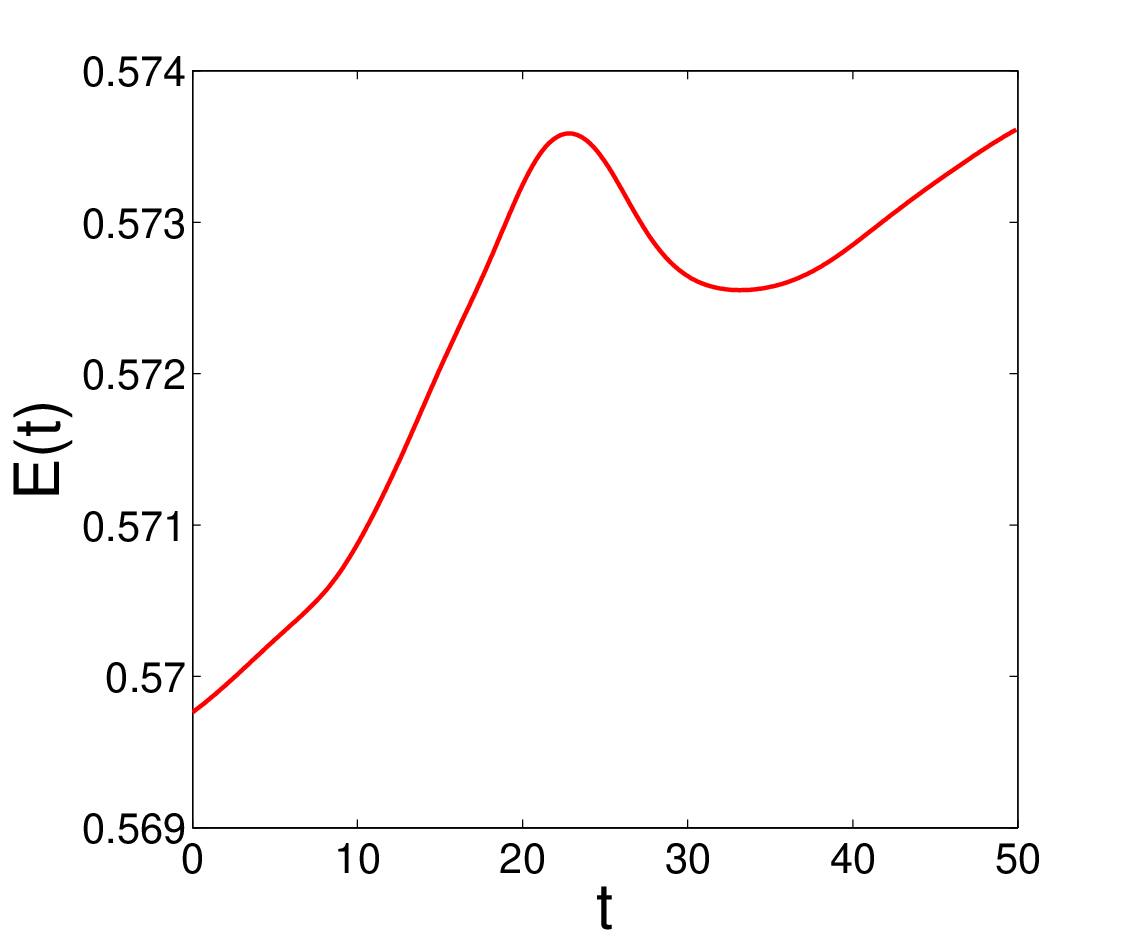}}
\subfigure{\includegraphics[width=1.6in,height=1.6in]{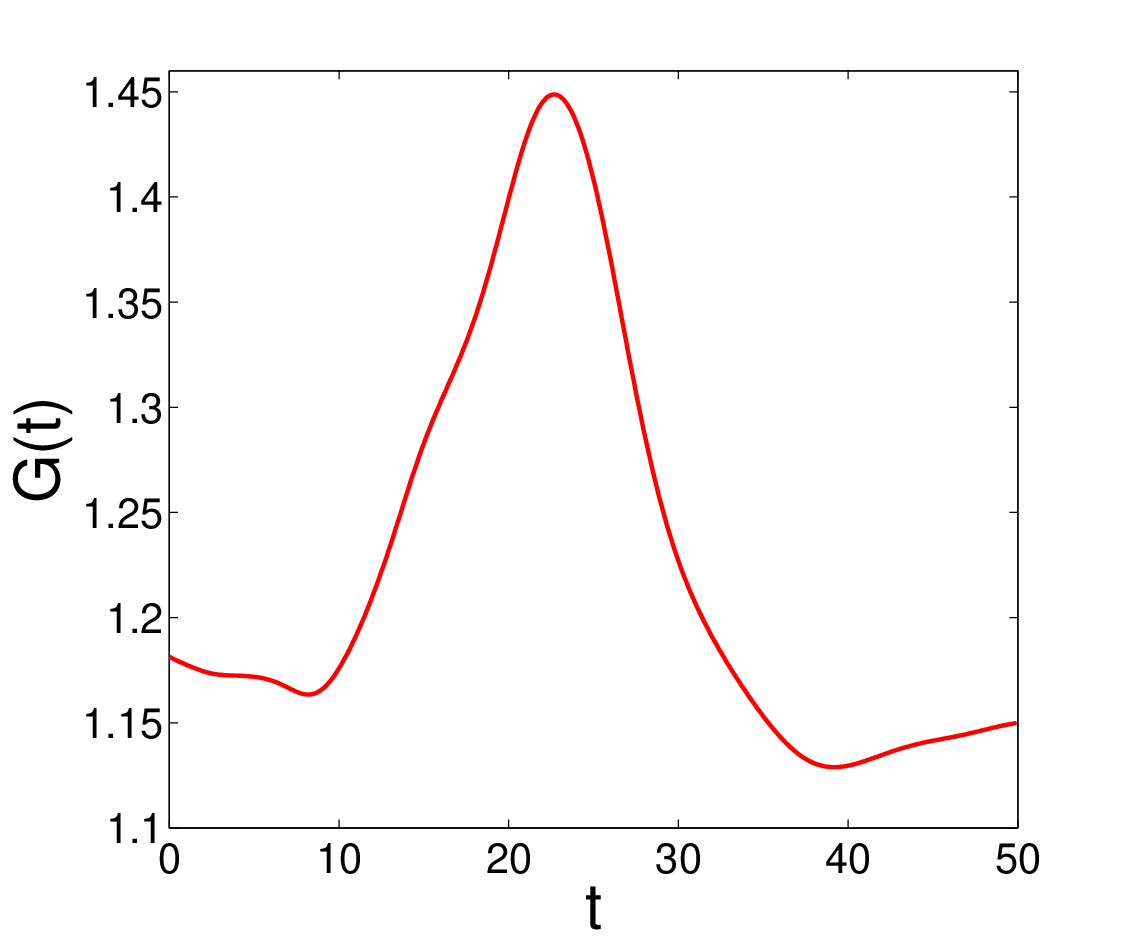}}
\subfigure{\includegraphics[width=1.6in,height=1.6in]{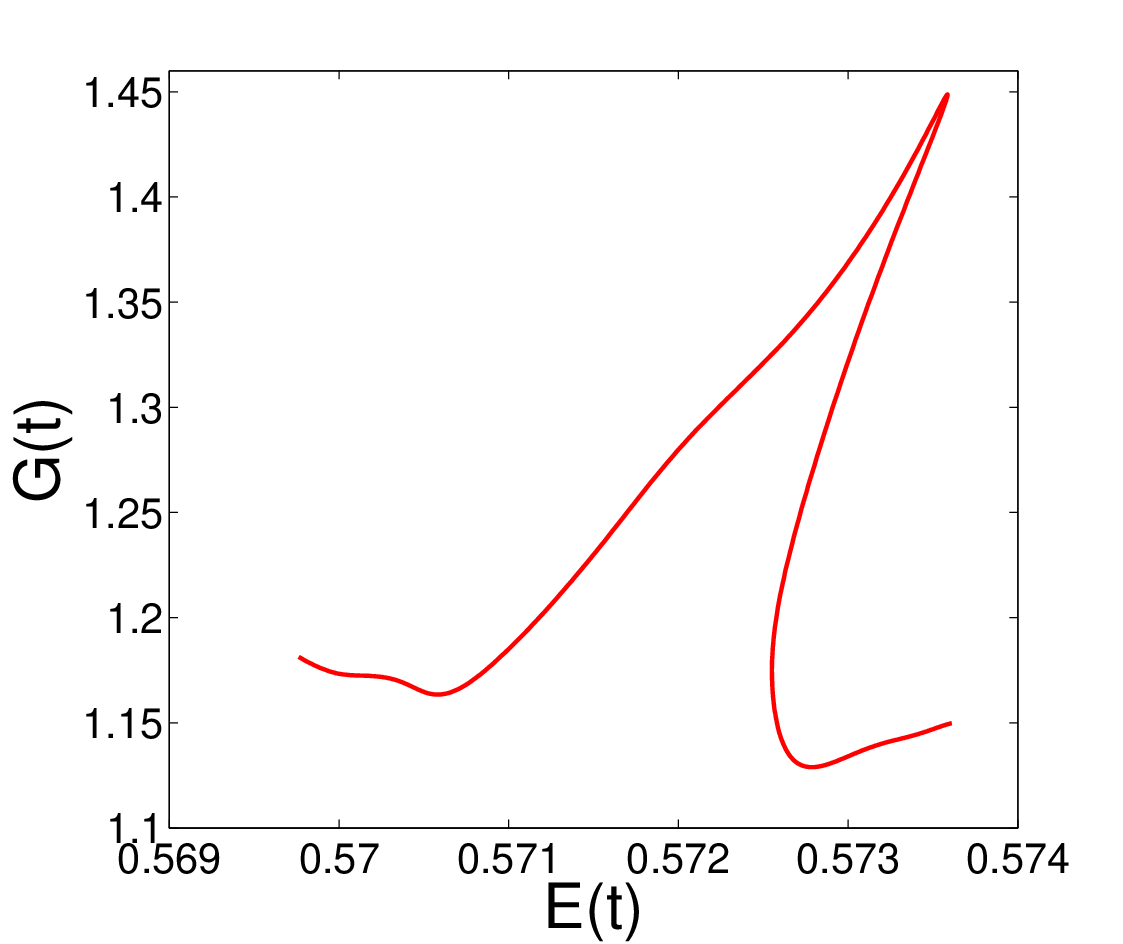}}
\subfigure{\includegraphics[width=1.6in,height=1.6in]{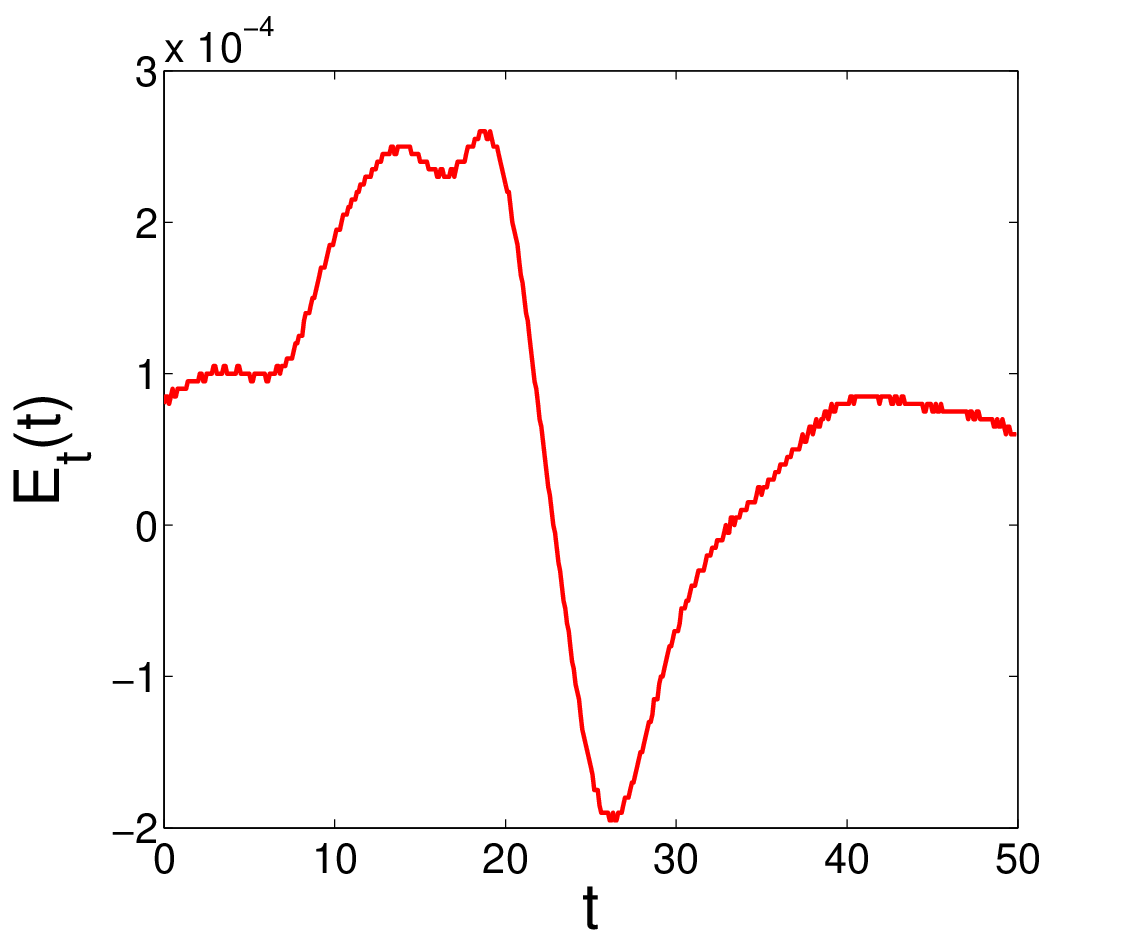}}
\subfigure{\includegraphics[width=1.6in,height=1.6in]{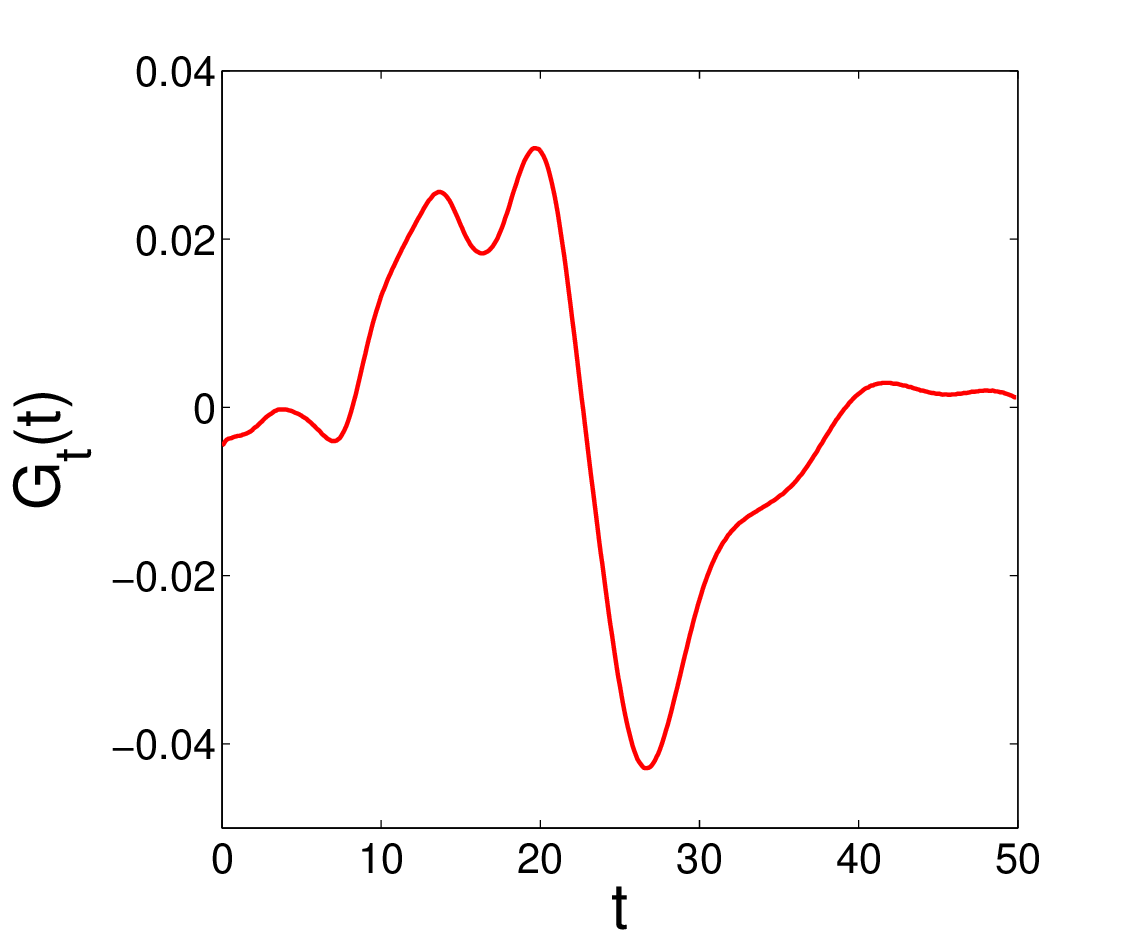}}
\subfigure{\includegraphics[width=1.6in,height=1.6in]{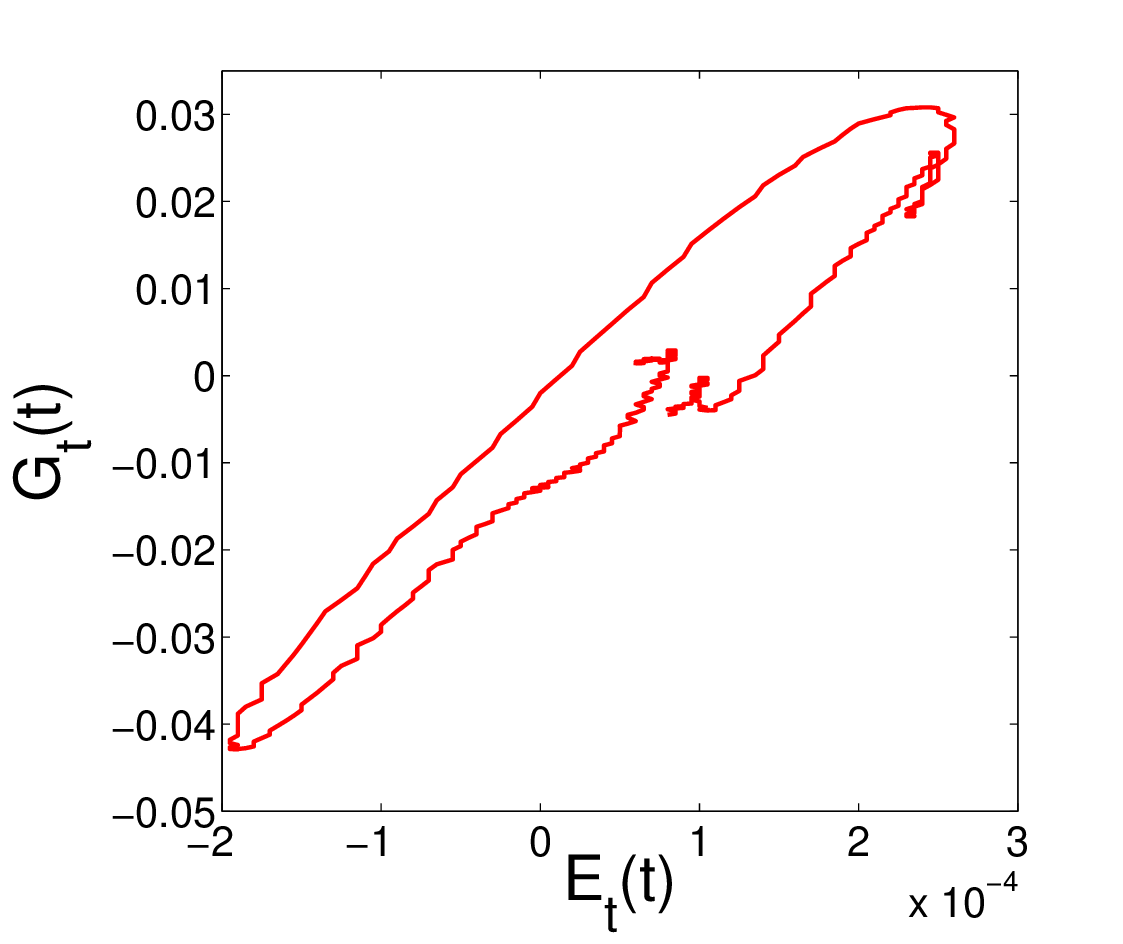}}
\caption{The modulation of kinetic energy $E(t)$ and enstrophy $G(t)$ in time $t$ when 
$n=1, R=10000$.}
\label{EG}
\end{figure}
\begin{figure}[ht] 
\centering
\subfigure[$t=0$]{\includegraphics[width=2.3in,height=2.3in]{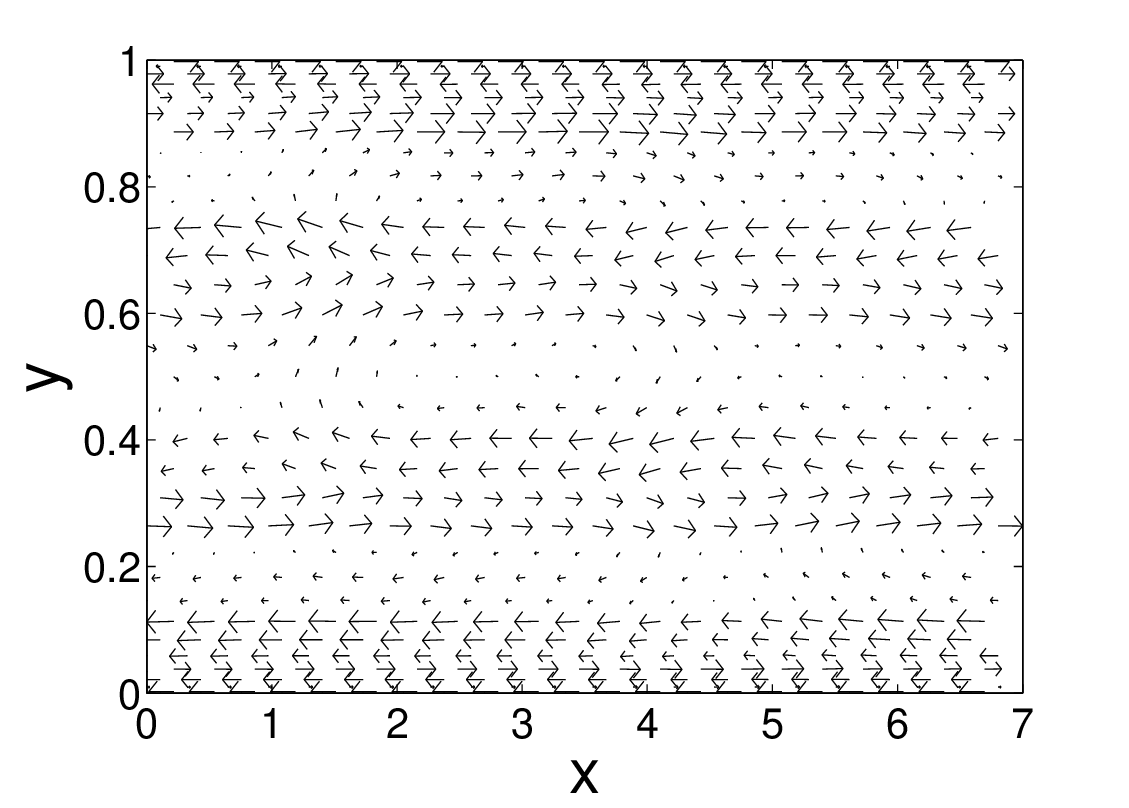}}
\subfigure[$t =4.9$]{\includegraphics[width=2.3in,height=2.3in]{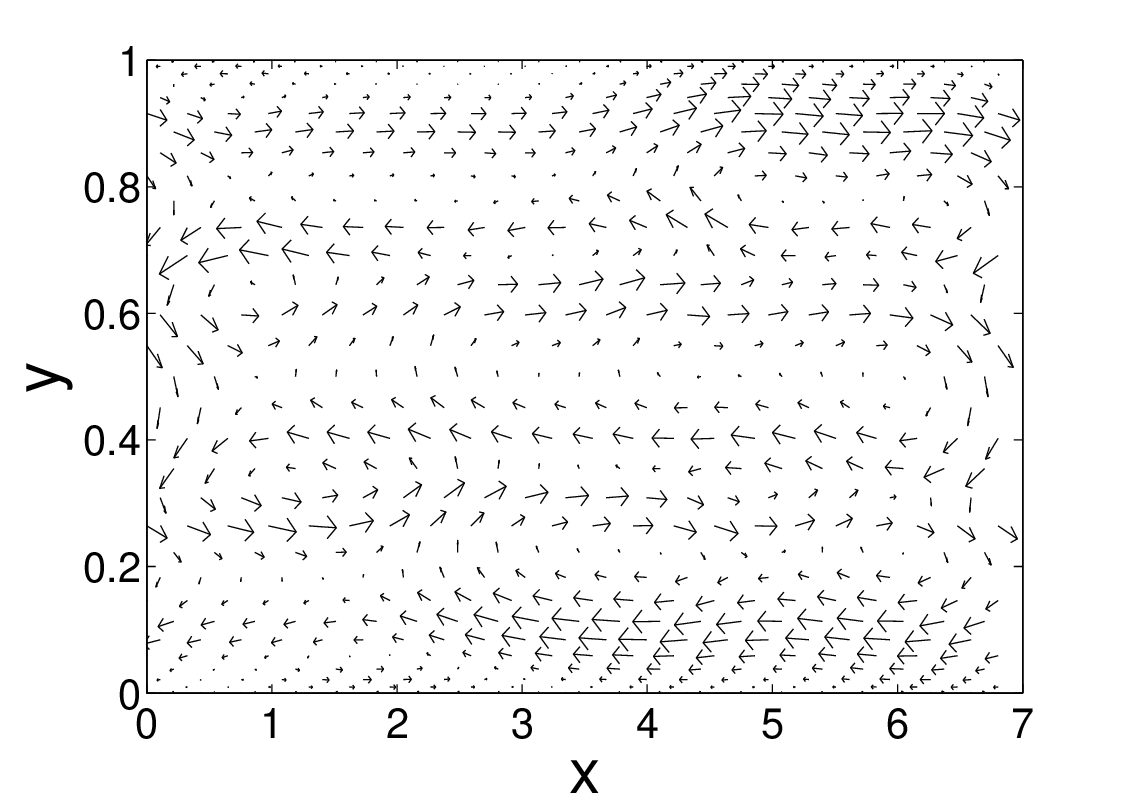}}
\subfigure[$t=24.9$]{\includegraphics[width=2.3in,height=2.3in]{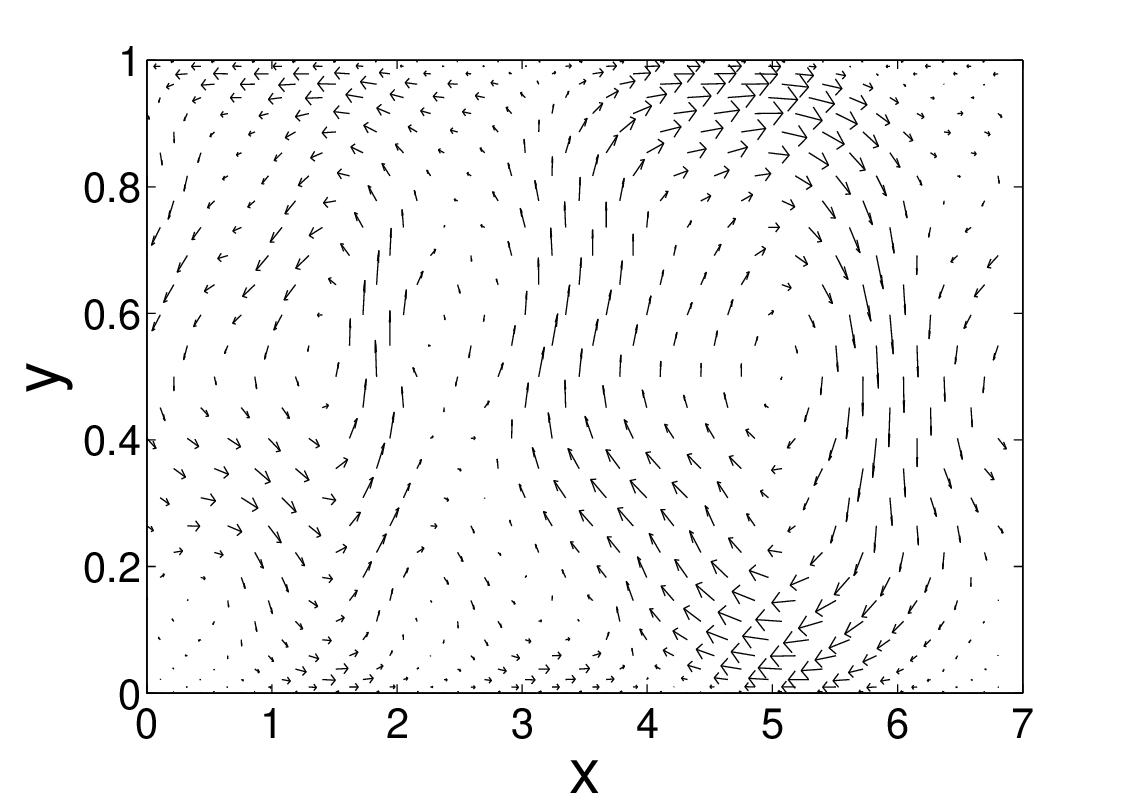}}
\subfigure[$L^2$ norm]{\includegraphics[width=2.3in,height=2.3in]{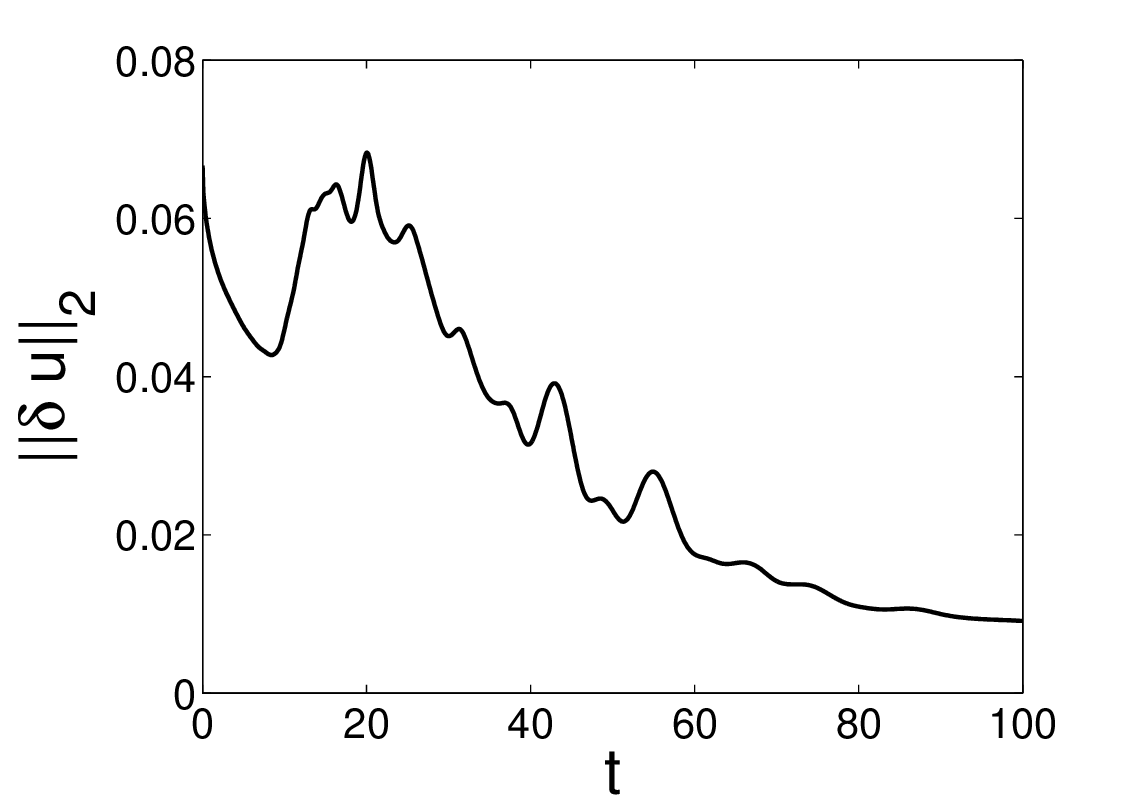}}
\caption{The development of a coherent structure and transient turbulence with initial condition (\ref{ICD}) and $R=10000$. Here the linear shear has been subtracted. (d) shows the growth of the $L^2$ norm of the deviation from the linear shear with initial condition (\ref{ICD}). The growth rate $\sg$ is defined as the $\log_e$ of quotient (of the first maximum and the first minimum) divided by the time spent. Here $\sg = 0.0293$.}
\label{DR}
\end{figure}
\begin{figure}[ht] 
\centering
\subfigure[$t=0$]{\includegraphics[width=2.3in,height=2.3in]{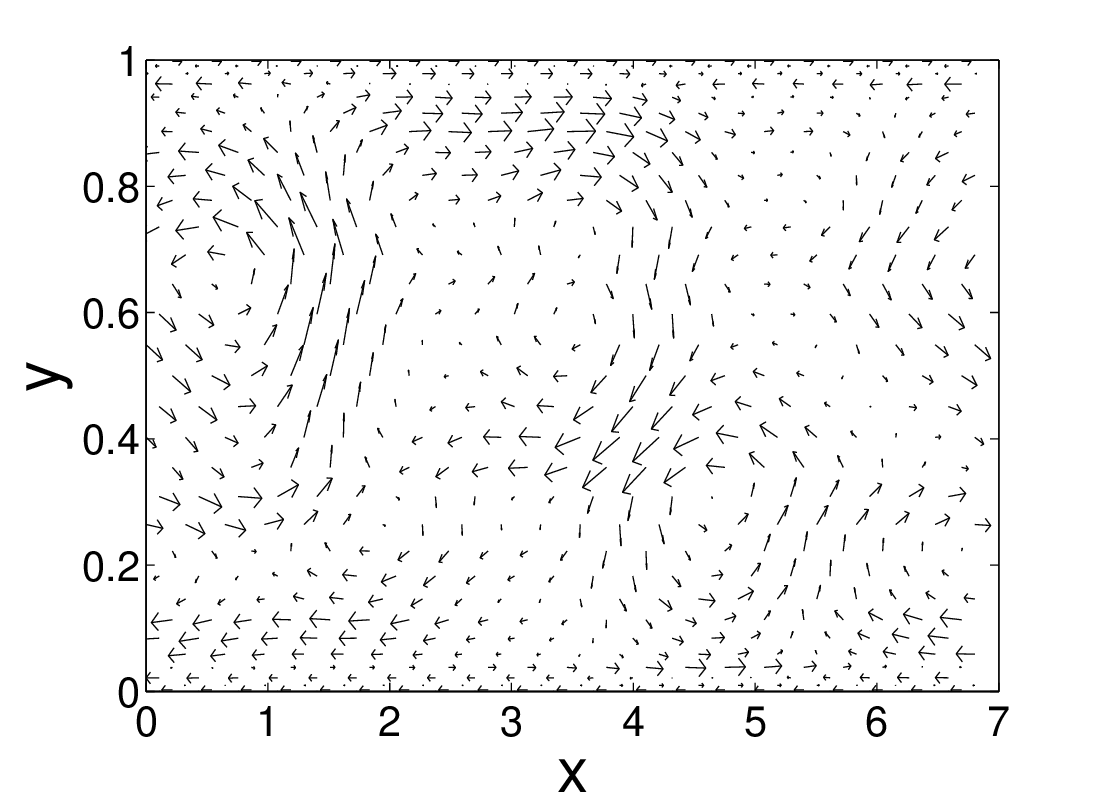}}
\subfigure[$t=4.9$]{\includegraphics[width=2.3in,height=2.3in]{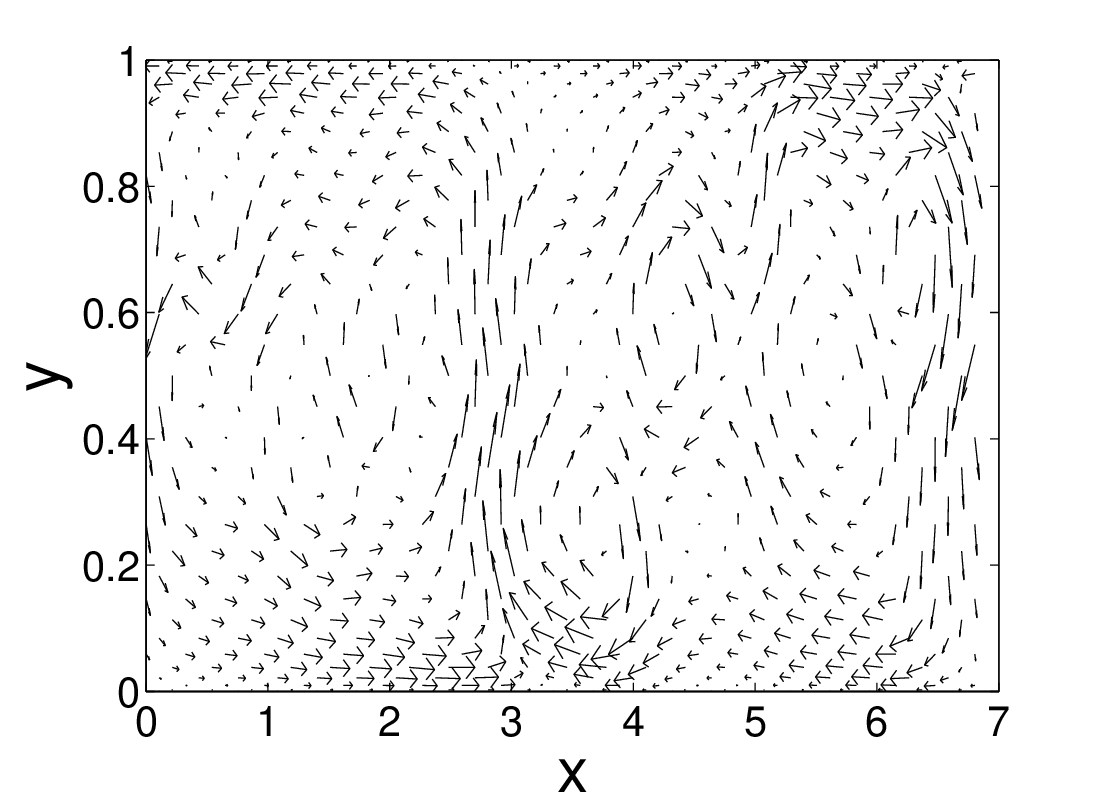}}
\subfigure[$t=14.9$]{\includegraphics[width=2.3in,height=2.3in]{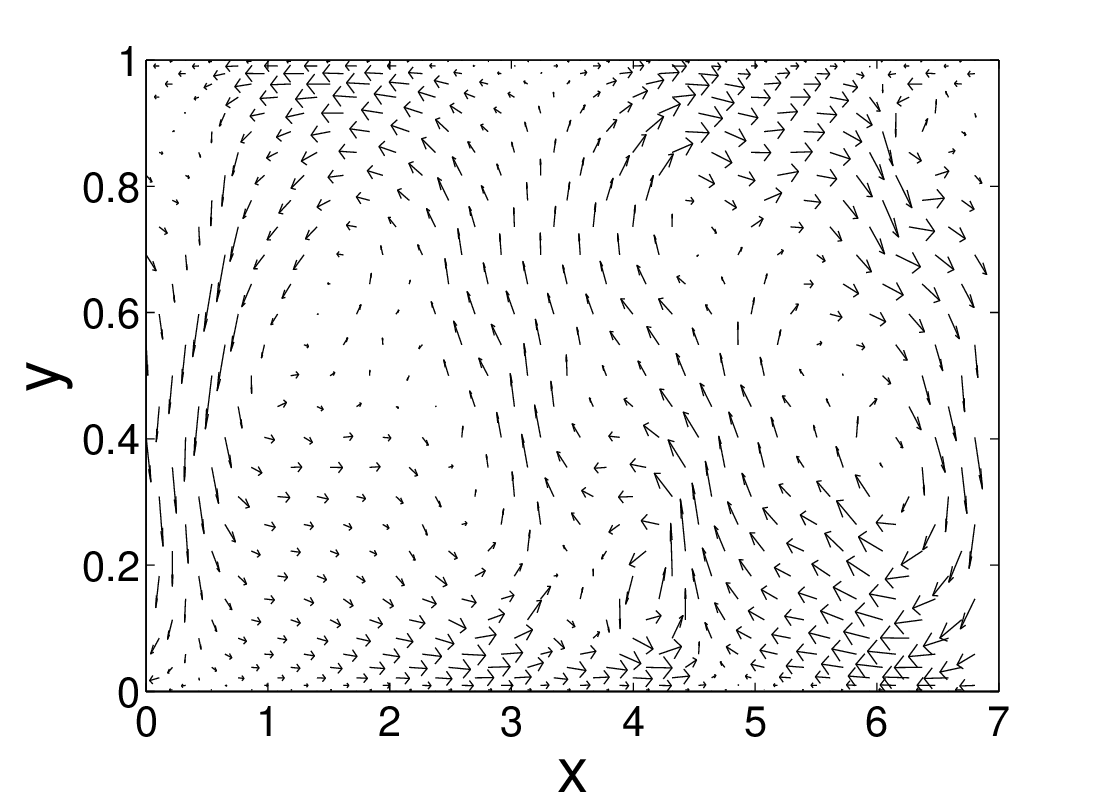}}
\subfigure[$t=24.9$]{\includegraphics[width=2.3in,height=2.3in]{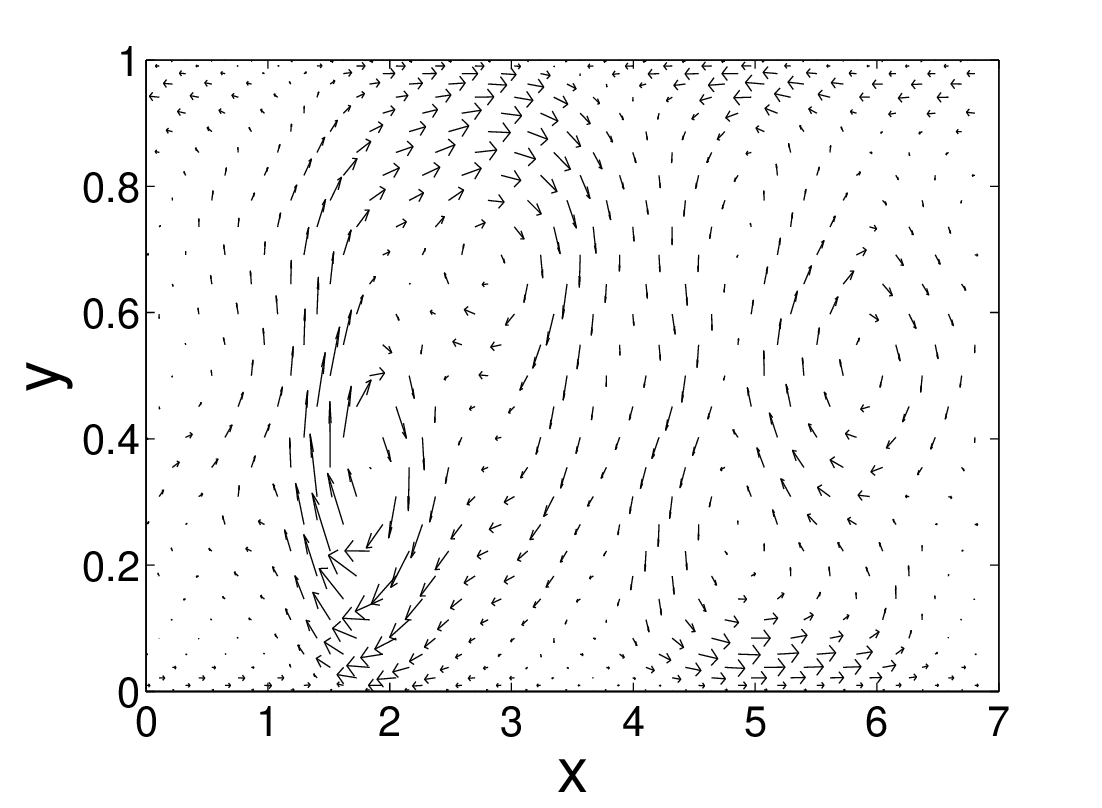}}
\caption{The development of a coherent structure and transient turbulence with initial condition (\ref{ICR}) and $R=10000$. Here the linear shear has been subtracted.}
\label{RF1}
\end{figure}

\begin{figure}[ht] 
\centering
\subfigure[$t=124.9$]{\includegraphics[width=2.3in,height=2.3in]{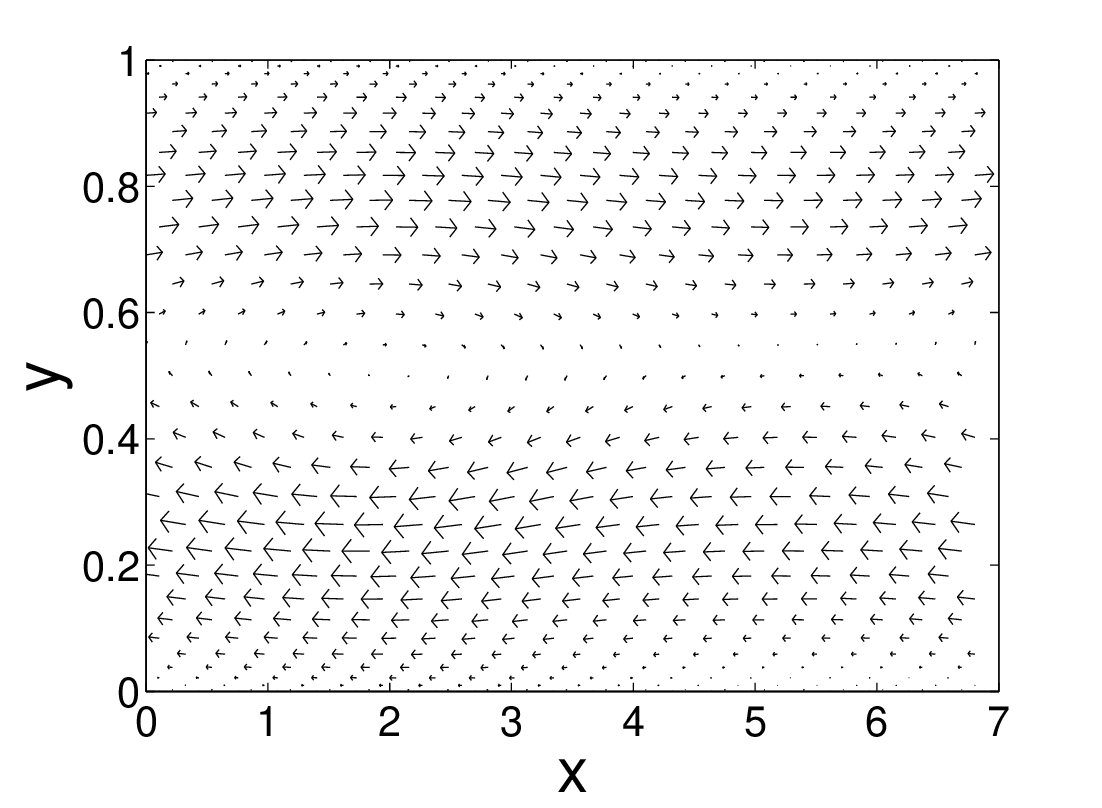}}
\subfigure[$L^2$ norm]{\includegraphics[width=2.3in,height=2.3in]{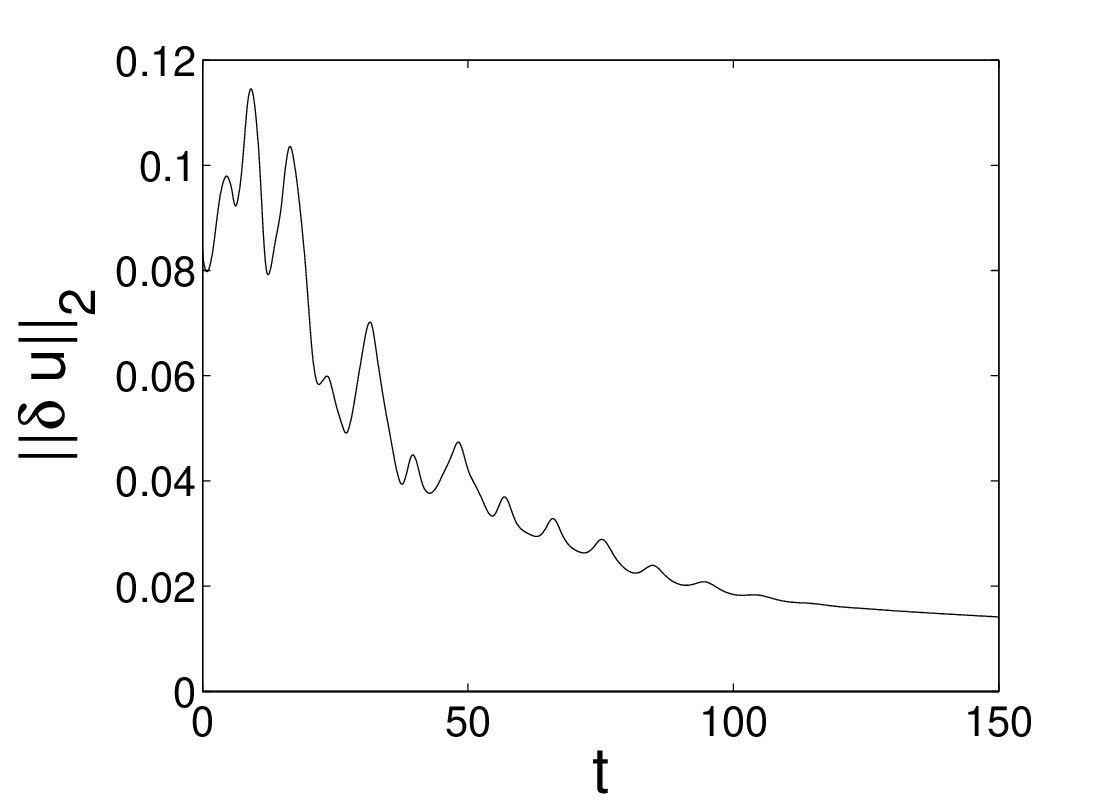}}
\caption{Figure \ref{RF1} continued. (b) shows the growth of the $L^2$ norm of the deviation from the linear shear with initial condition (\ref{ICR}). The growth rate $\sg$ is defined as the $\log_e$ of quotient (of the first maximum and the first minimum) divided by the time spent. Here $\sg = 0.04$.}
\label{RF2}
\end{figure}

\begin{figure}[ht] 
\centering
\subfigure[$t=0$]{\includegraphics[width=2.3in,height=2.3in]{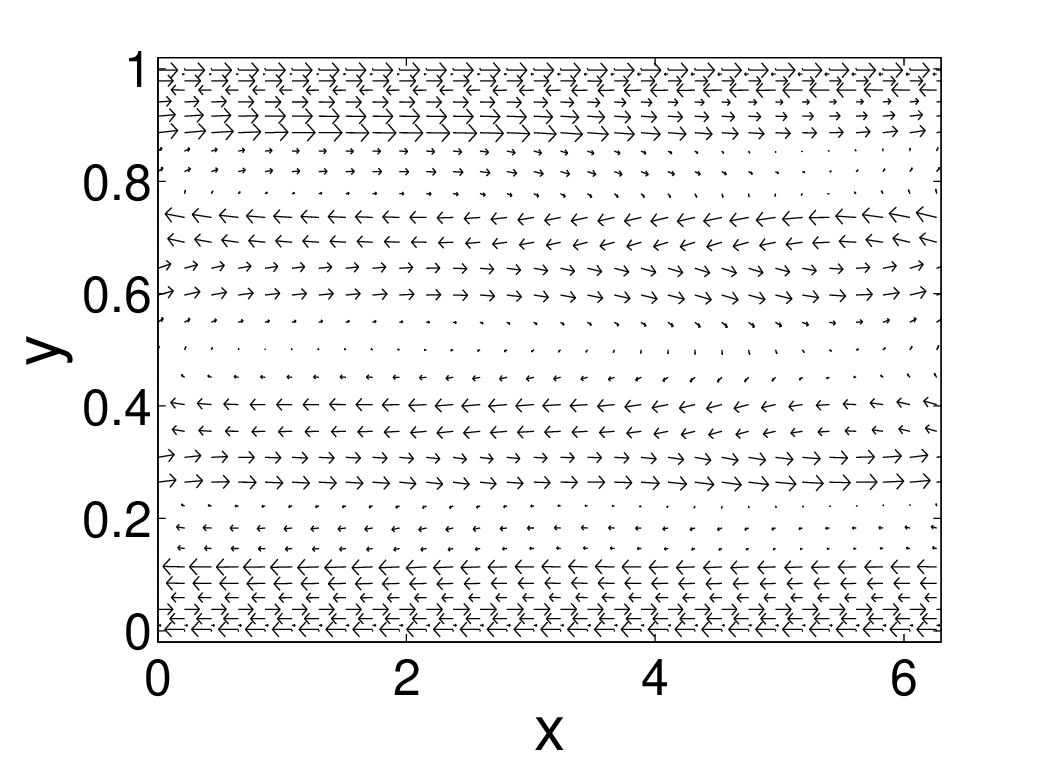}}
\subfigure[$t=12.5$]{\includegraphics[width=2.3in,height=2.3in]{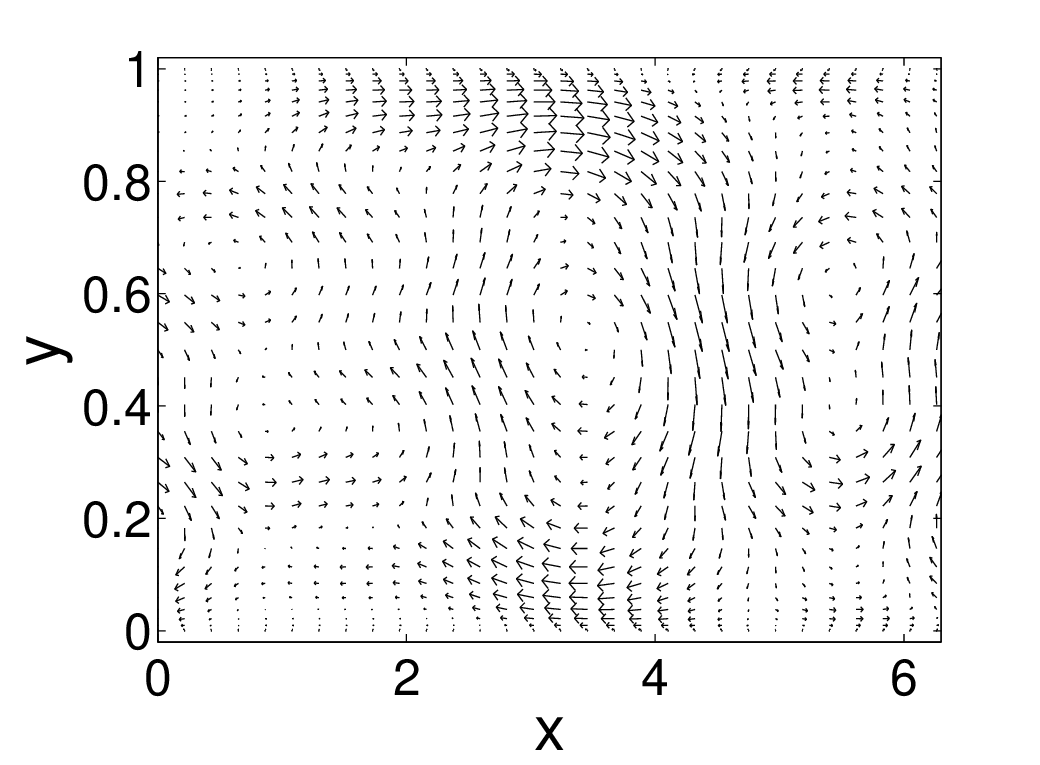}}
\subfigure[$t=25$]{\includegraphics[width=2.3in,height=2.3in]{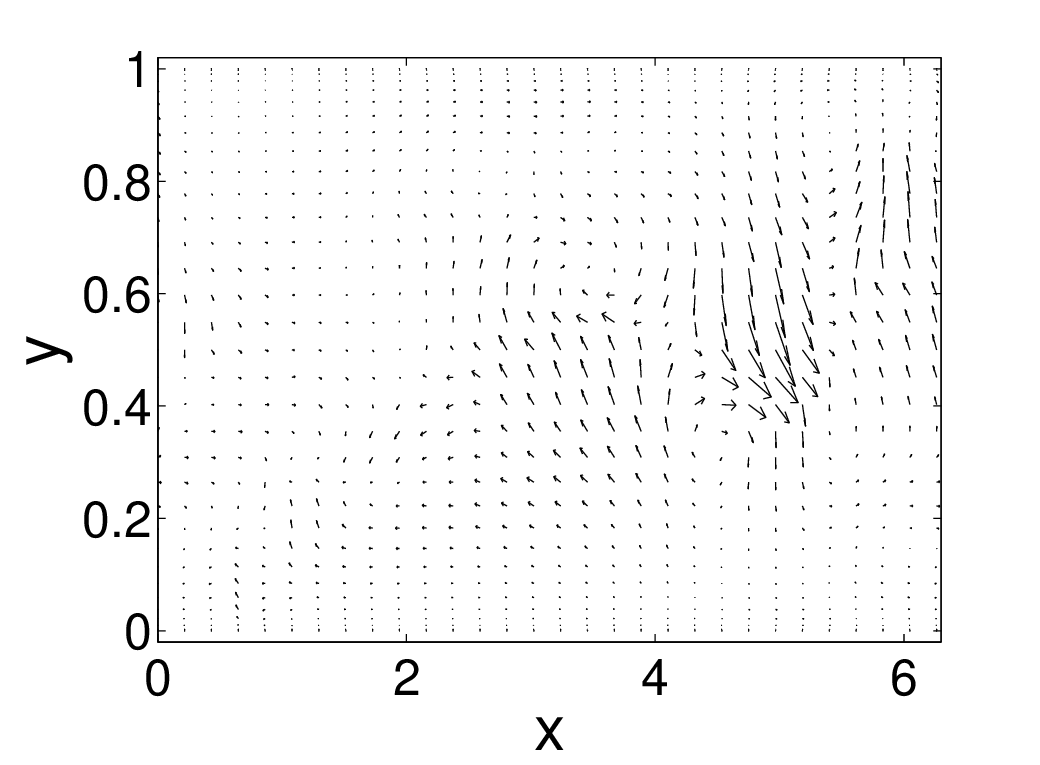}}
\subfigure[$t=35$]{\includegraphics[width=2.3in,height=2.3in]{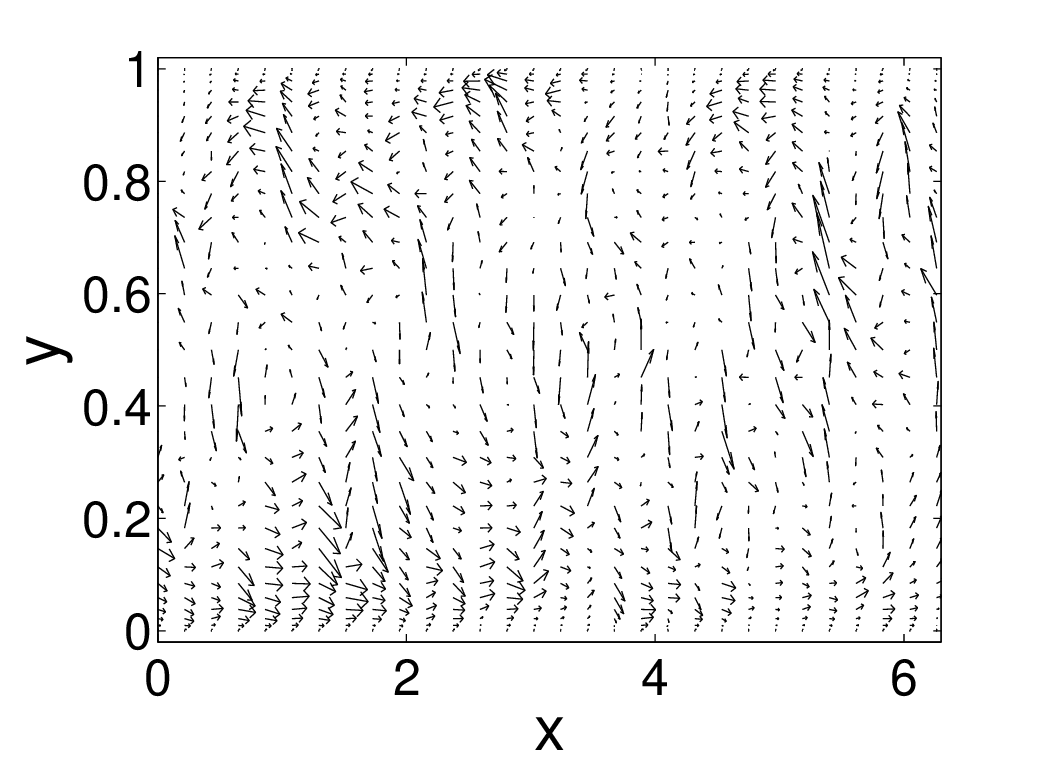}}
\caption{The development of a coherent structure and transient turbulence in the ($x,y$)-section at $z=0.6\pi$, 
with initial condition (\ref{3ICD}) and $R=5000$. Here the linear shear has been subtracted.}
\label{3DF1}
\end{figure}

\begin{figure}[ht] 
\centering
\subfigure[$t=150$]{\includegraphics[width=2.3in,height=2.3in]{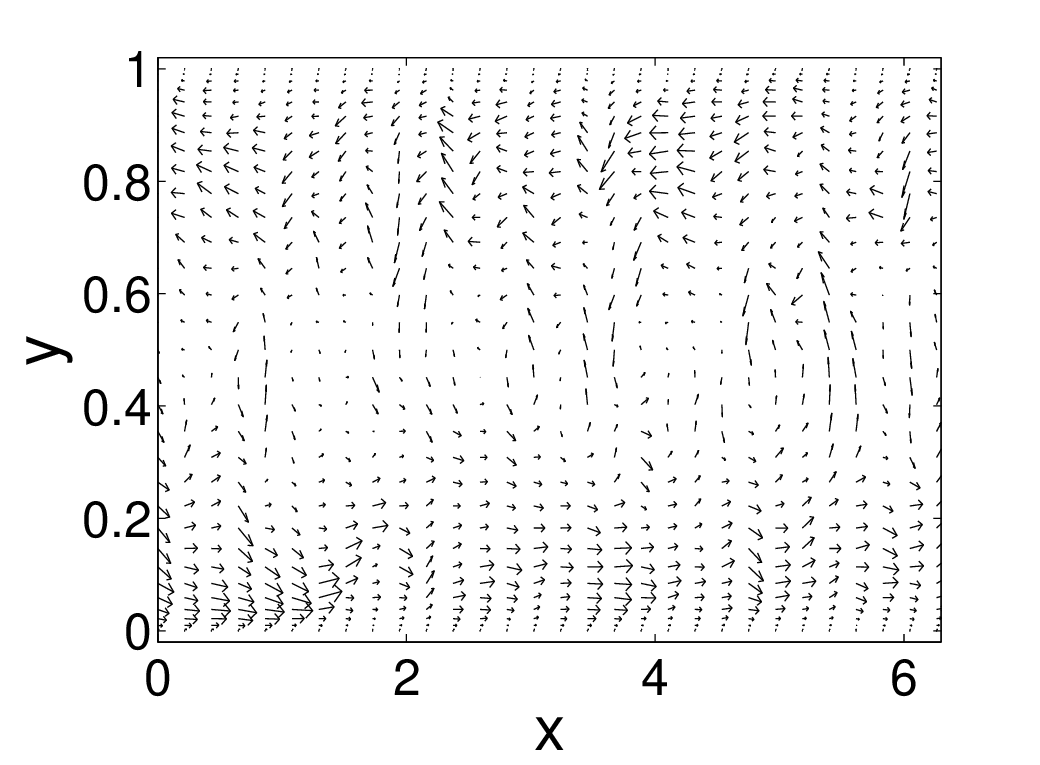}}
\subfigure[$L^2$ norm]{\includegraphics[width=2.3in,height=2.3in]{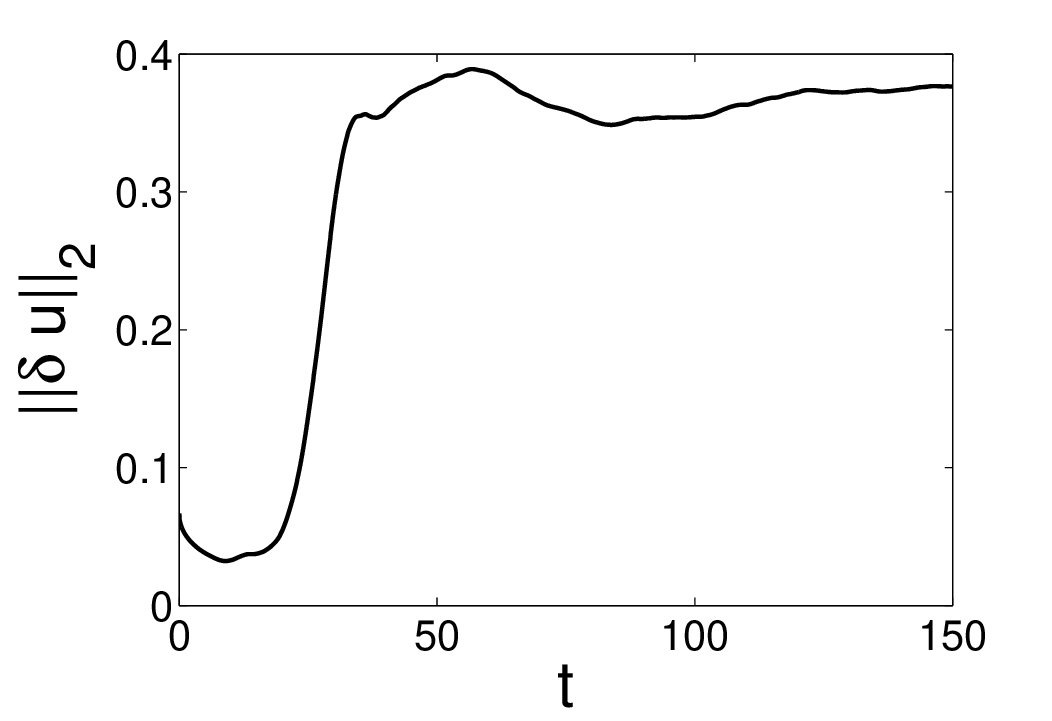}}
\caption{Figure \ref{3DF1} continued. (b) shows the growth of the $L^2$ norm of the deviation from the linear shear with initial condition (\ref{3ICD}).}
\label{3DF2}
\end{figure}

\begin{figure}[ht] 
\centering
\subfigure[$t=0$]{\includegraphics[width=2.3in,height=2.3in]{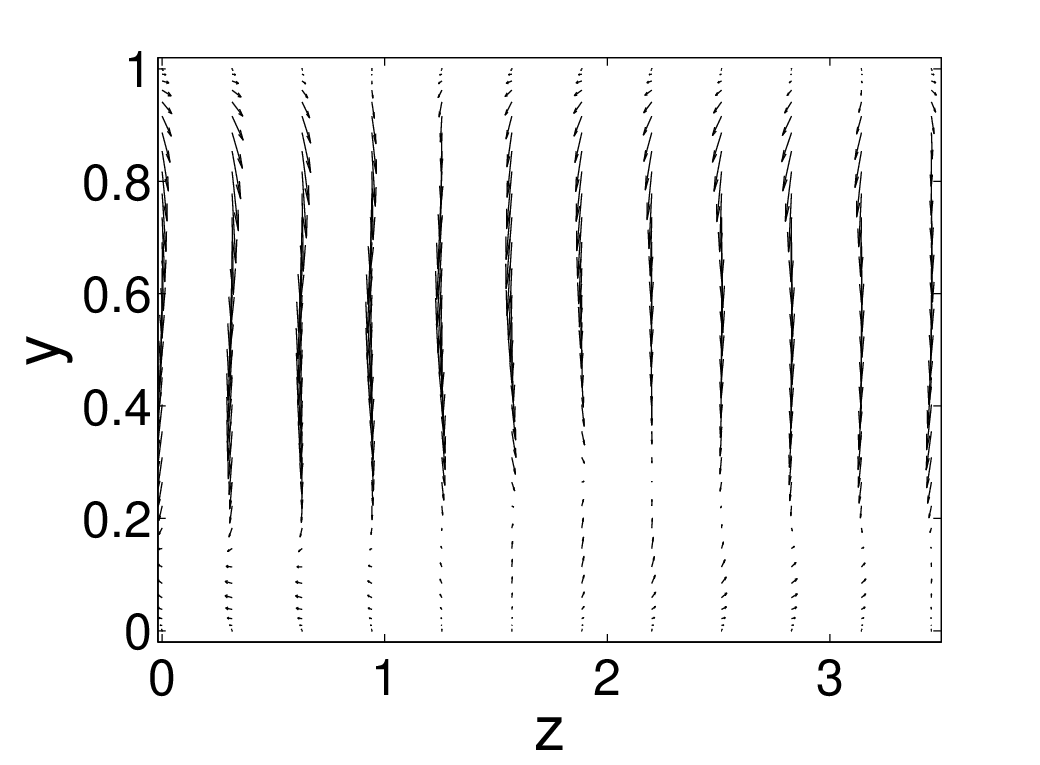}}
\subfigure[$t=12.5$]{\includegraphics[width=2.3in,height=2.3in]{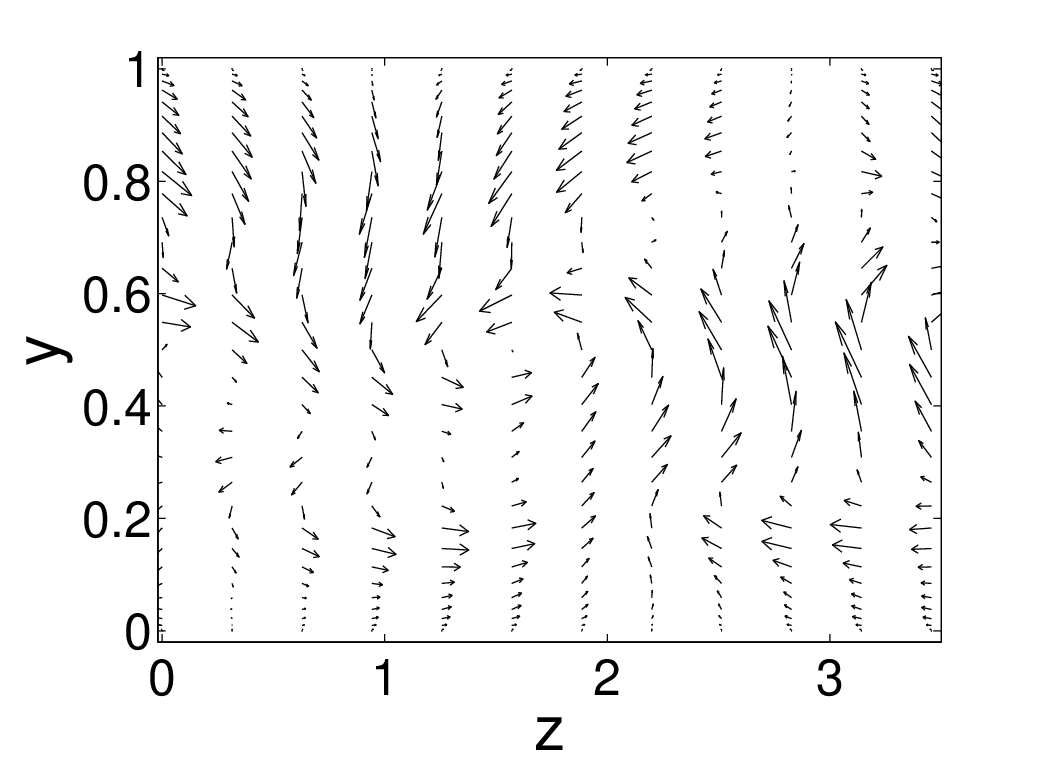}}
\subfigure[$t=35$]{\includegraphics[width=2.3in,height=2.3in]{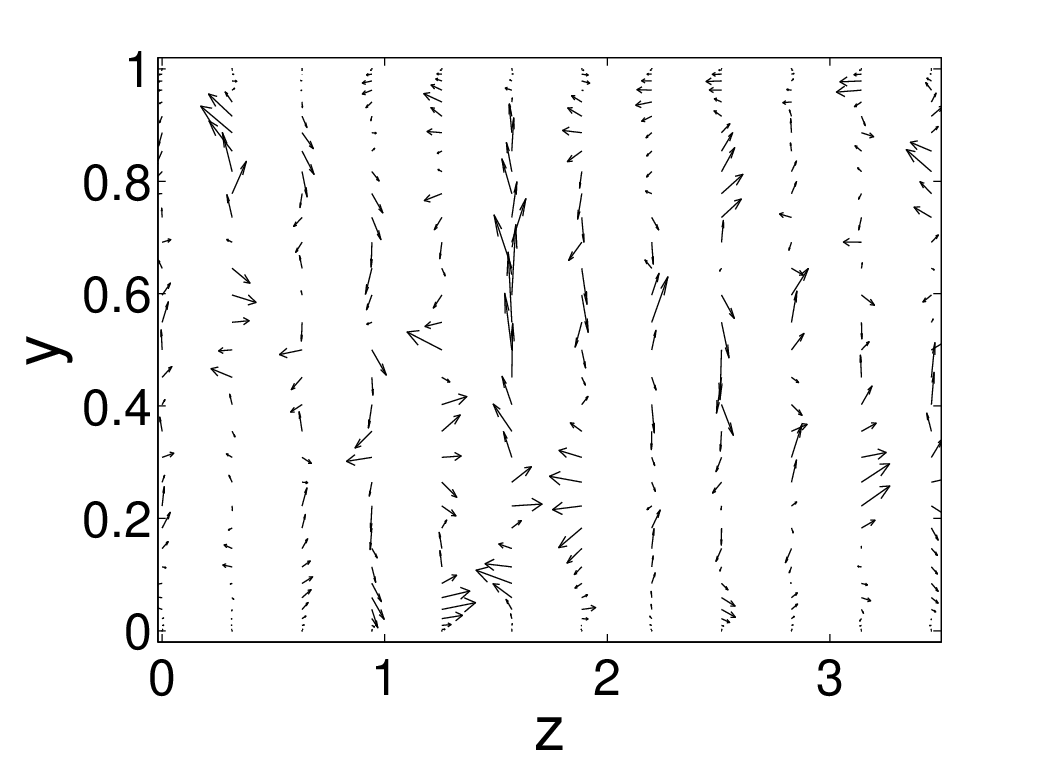}}
\subfigure[$t=150$]{\includegraphics[width=2.3in,height=2.3in]{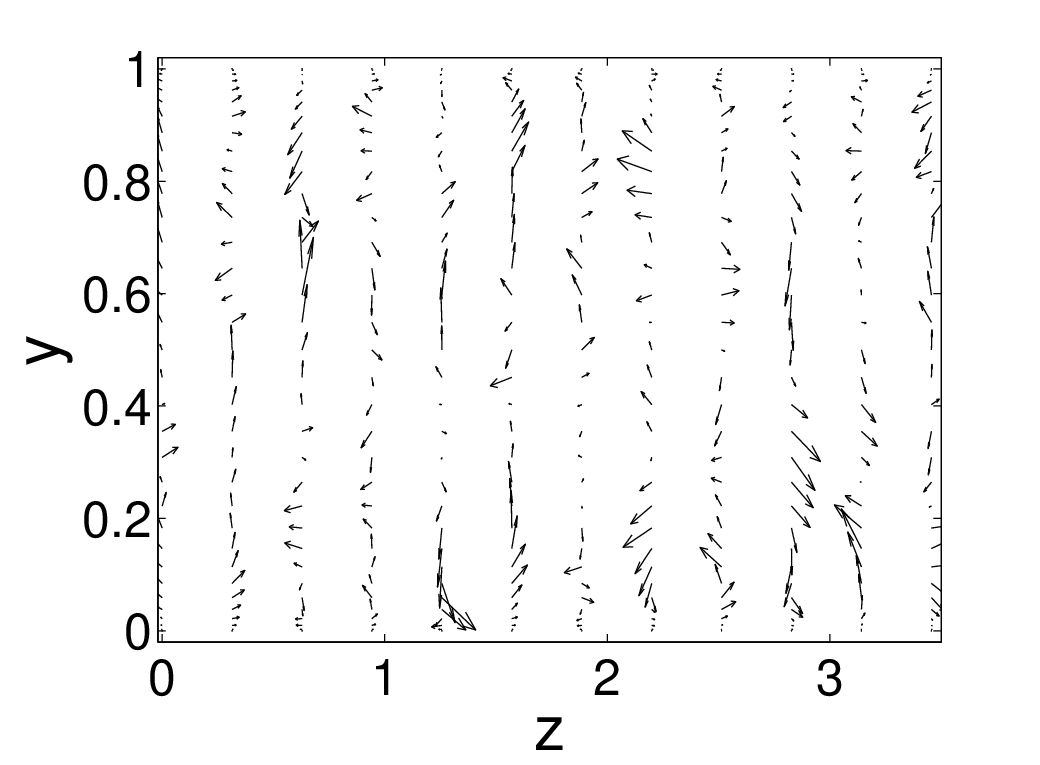}}
\caption{The development of a coherent structure and transient turbulence in the ($y,z$)-section at $x=\frac{34.1}{32}\pi$, 
with initial condition (\ref{3ICD}) and $R=5000$.}
\label{3DF3}
\end{figure}

In Figure \ref{cod}, we show a few time flashes of the velocity field evolution. One can see
clearly the development of a coherent structure which is the viscous continuation of the inviscid cat's eye steady 
state shown in \cite{LL09}. In Figure \ref{grth}, we show the velocity 
$L^2$ norm deviation of the solution from the linear shear. One can see the development of pulses. We believe that 
the center-unstable and center-stable manifolds (if they exist) of the 2D Euler equations coincide. The slightly 
viscous continuations of them deviate slightly. By the mechanism of $\la$-lemma \cite{PM82} \cite{Li03a}, generic 
orbits in the neighborhood of the drifting shear (\ref{sdr}) approach the 
center-unstable manifold, and this corresponds to the time duration before the peak of the first pulse. The orbit 
first follows the center-unstable manifold. Since the center-unstable and center-stable manifolds deviate slightly, 
after reaching the first peak, the orbit then follows the center-stable manifold to return. Repeating the same process, 
more pulses may be generated. 
Comparing Figure \ref{cod} with Figure \ref{grth}(a), one can see that the development of the coherent structure in 
Figure \ref{cod} corresponds to the development of the first pulse in Figure \ref{grth}(a). We can measure the growth rate 
of the uphill of the first pulse. For example in Figure \ref{grth}(a), let $m$ be the first minimum, $M$ be the first 
maximum, and $\Dl t$ be the duration from the first minimum to the first maximum, then the exponential growth rate
$\sg$ is defined by
\begin{equation}
\sg = \frac{1}{\Dl t} \ln \frac{M}{m} . 
\label{GR}
\end{equation}
This exponential growth rate measures the exponential growth in time of the deviation from the linear shear.
In the cases we simulated, the exponential growth rate can be as large as $\sg = 0.25$. The exponential growth rate 
here can also be regarded as a measure on how fast the flow is initially leaving the linear shear to transient 
turbulence.
In Figure \ref{grth2}, we show the velocity $L^2$ norm deviation of the same solution but from the moving frame of 
the slow drifting of the oscillatory shear (\ref{sdr}). We denote by $u_0$ the velocity field given by (\ref{sdr}). 
Comparing Figure \ref{grth} and Figure \ref{grth2}, one can see that the time instant of the peak of the first 
pulse agrees quite well for $n=1,2$. For $n=3$, the initial random perturbations almost overwhelm the oscillatory 
part of the oscillatory shear. That is why the agreement is not that good. Such an agreement shows that the first 
pulse is generated by the linear instability of our oscillatory shear; and the slow drifting (\ref{sdr}) has very 
little effect on the development of the first pulse. 
From (\ref{cr}) and (\ref{sdr}), we see that the time $T$ for the oscillatory shear to drift outside its unstable regime 
(\ref{cr}) can be estimated by 
\[
ce^{-\e (4n\pi )^2 T} = \frac{1}{2}\frac{1}{4\pi}.
\]
That is,
\begin{equation}
T = \frac{1}{\e (4n\pi )^2}\ln \frac{7}{4} .
\label{TS}
\end{equation}
This time scale agrees well with the time scale of the first pulse in the $L^2$ norm evolution of the random 
perturbations. For example, in the setting of Figure \ref{grth2}(a), $T\approx 35$ which almost coincides with the 
duration of the first pulse. That is why the flow decays after the first pulse. For Figure \ref{grth}(b), $T \approx 70$
which coincides with the duration of the 3 pulses. 
We can also measure the exponential growth rate of the uphill 
of the first pulse as for Figure \ref{grth} using the same definition (\ref{GR}). Here the exponential growth rate 
approximates the unstable eigenvalue of our oscillatory shear. One can see that all the exponential growth rates 
in Figure \ref{grth2} are very close to each other. When $n$ is the same and the Reynolds number $R$ is different, 
this fact is expected since the unstable eigenvalue of our oscillatory shear approaches its inviscid unstable 
eigenvalue as $R \ra \infty$. The exponential 
growth rate here measures how fast the initial random perturbation grows when observed on the moving frame of 
the slow drifting (\ref{sdr}) toward the linear shear. The uphill of the pulse is generated by nonlinear growth induced by the unstable 
eigenvalues. Since our Reynolds number is quite large, the Navier-Stokes equation is ``near" 
the Euler equation, and there is a stable eigenvalue corresponding to an unstable eigenvalue 
due to the near conservativeness. The stable eigenvalue corresponds to the downhill of 
the pulse. Concerning other pulses in Figures \ref{grth} and \ref{grth2}, e.g. Figure \ref{grth2}(e)(f), since the random 
perturbations in these cases almost overwhelm the oscillatory part of our oscillatory shear, other instability 
sources play significant roles. The first pulse is due to the linear instability of our oscillatory shear. 
When the first pulse is finished, our oscillatory shear has drifted outside its unstable regime (\ref{sdr}), but 
the new transient state can pick further instability and develop even higher pulses. This further illustrates our point 
that the linear instability of our oscillatory shear serves as the initiator for transition to turbulence. 

Finally, it is always tempting to try to apply the tools of dynamical systems to characterize turbulence. One promising 
tool is the Melnikov integral. In \cite{LL08}, we built the Melnikov integral for the 2D Kolmogorov flow (i.e. with 
periodic boundary conditions in both spatial directions and an artificial force). The Melnikov integral there was built from the kinetic energy and enstrophy, and evaluated along a heteroclinic orbit. For the Couette flow, we have not found any proper heteroclinic orbit. Nevertheless, we can still investigate the modulation of the kinetic energy and enstrophy. When $\e = 0$, the boundary condition (\ref{BC}) reduces to just the non-penetrating condition
\[
u_2(t, x, 0) = u_2(t, x, 1) = 0 . 
\]
A direct calculation shows that the kinetic energy
\[
E=\int (u_1^2 +u_2^2) dx dy
\]
and the enstrophy
\[
G=\int \Om^2 dxdy , \quad \quad \text{ where } \Om = \pa_x u_2 - \pa_y u_1 
\]
are invariant under the Euler dynamics. When $\e >0$, they are no longer invariant, but their time derivatives
$E_t$ and $G_t$ should be small for small $\e$: 
\[
E_t = 2 \e  \int  u_i u_{i,jj} dxdy , \quad G_t = 2 \e  \int   \Om \Om_{,jj} dxdy .
\]
In Figure \ref{EG}, we plot the evolution of $E$, $G$, $E_t$ and $G_t$ along the orbit in Figure \ref{cod}. Due to 
the high-frequency oscillation nature of the oscillatory part of our oscillatory shear (\ref{Os}), the modulation 
of the enstrophy (characterized by $G_t$) is about $100$ times of that of the kinetic energy. This shows that 
during the initiation of transition to turbulence, there has been a massive vorticity production.

\section{Other numerical results}

\subsection{Random initial perturbations with dominant random shears}

We choose the x-direction spatial period $L_x = 2.2 \pi$, and the initial condition
\begin{equation}
u_1 = y + \tU (y) + \tu_1(x,y),\quad u_2 = \tu_2 (x,y);
\label{ICD}
\end{equation}
where $\tU$, $\tu_1$ and $\tu_2$ are all random perturbations; the 
$L^2$ norm of $\tU$ is $0.1$, and the $L^2$ norms of $\tu_1$ and $\tu_2$ are
$0.01$. We choose the Reynolds number $R=10000$. In comparison with Figure \ref{cod}, a 
similar coherent structure is developed, see Figure \ref{DR}. This shows that the transition 
to turbulence starting from the initial condition (\ref{ICD}) is similar to the transition from the 
initial condition in Figure \ref{cod}. The order $0.1$ term of the initial condition (\ref{ICD}) is 
a 2D shear which should have the similar linear instability as the sequence of oscillatory shears (\ref{Os}).

\subsection{Random initial perturbations without dominant random shears}

We choose the x-direction spatial period $L_x = 2.2 \pi$, and the initial condition
\begin{equation}
u_1 = y + \tu_1(x,y),\quad u_2 = \tu_2 (x,y);
\label{ICR}
\end{equation}
where $\tu_1$ and $\tu_2$ are all random perturbations, and the 
$L^2$ norms of $\tu_1$ and $\tu_2$ are $0.1$. We choose the Reynolds number $R=10000$.
Even though the initial condition (\ref{ICR}) has no dominant shear, a 
similar coherent structure as those of Figures \ref{cod} and \ref{DR} is developed, see Figures \ref{RF1} 
and \ref{RF2}. At $t=0$, the velocity field in Figure \ref{RF1} is very different from those of 
Figures \ref{cod} and \ref{DR}. On the other hand, at $t=14.9$, Figure \ref{RF1}(c) and Figure \ref{cod}(c) are 
very similar. It is also interesting to notice that the velocity field starting from (\ref{ICR}) decays to 
a near 2D shear state later on, see Figure \ref{RF2}(a). This numerical experiment (Figures \ref{RF1} and \ref{RF2}) 
further illustrates the role of linearly unstable slow orbits in the neighborhood of the linear shear, as the 
initiators for the transition to turbulence.

\subsection{3D random initial perturbations with dominant random shears}

We choose the x-direction spatial period $L_x = 2.2 \pi$, z-direction spatial period $L_z = 1.2 \pi$, the Reynolds number
$R=5000$, and the initial condition

\begin{equation}
u_1 = y + \tU (y) + \tu_1(x,y,z),\quad u_2 = \tu_2 (x,y,z),\quad u_3 = \tu_3 (x,y,z);
\label{3ICD}
\end{equation}
where $\tU$, $\tu_1$, $\tu_2$ and $\tu_3$ are all random perturbations; the 
$L^2$ norm of $\tU$ is $0.1$, and the $L^2$ norms of $\tu_1$, $\tu_2$ and $\tu_3$ are
$0.01$. In Figures \ref{3DF1} and \ref{3DF2}, we show the velocity field development of 
the ($x,y$)-section at $z=0.6\pi$. We observe that a similar coherent structure as in the 2D case before, is developed 
during the initial stage of the development, e.g. $t\in [0,25)$. This time interval corresponds to the valley and the first 
half of the uphill slope in Figure \ref{3DF2}(b). During the second half of the uphill slope, the velocity field transits into 
more violent and permanent turbulent state. The time $t=35$ (Figure \ref{3DF1}(d)) is near the top of the uphill slope. At $t=35$, 
the velocity field enters into more fully developed turbulence. We believe that, unlike 2D Euler equations, center-unstable and
 center-stable manifolds (if they exist) of the 3D Euler equations do not coincide due to the non-conservation of their enstrophy. 
Thus, under 3D slightly viscous dynamics, orbits following the center-unstable manifold by the $\la$-lemma, will not follow the 
center-stable manifold to return, rather follow the center-unstable manifold into fully developed turbulence, as shown in 
Figure \ref{3DF2} for $t\in [35,150]$. In fact, the center-unstable manifold is far away from the center-stable manifold, and 
the center-unstable manifold connects to the fully developed turbulence network. In Figure \ref{3DF3}, we show the corresponding 
velocity field development of the ($y,z$)-section at $x = \frac{34.1}{32}\pi$. Here the vector fields are rescaled, e.g. the true 
scale of the vector field in Figure \ref{3DF3}(a) is of order $0.01$. In summary, this numerial experiment with initial 
condition (\ref{3ICD}) indicates that the linear instability of the dominant 2D shear initiates the transition to 3D turbulence. In 
3D, slow orbits in the neighborhood of the linear shear are more abundant. In particular, there are two types of 3D shears which 
can be more unstable than 2D shears \cite{Li11a} \cite{Li11b}. Thus in 3D, there are more initiators for the transition from the 
linear shear to turbulence.

\section{Conclusion}

It is claimed in \cite{LL09} on the resolution of the Sommerfeld paradox, that linearly unstable slow orbits in arbitrarily 
small neighborhoods of the Couette linear shear are the initiators for the transition to turbulence, and their linear 
instabilities should generate transient nonlinear growth. The current numerical simulations verified this for both 2D and 3D 
transitions.

\end{document}